\documentclass[11pt,a4paper]{article}
\usepackage{amssymb,amsmath,amsfonts,mathrsfs,bm,enumerate,comment,overpic}
\allowdisplaybreaks[1]
\numberwithin{equation}{section}

\usepackage[colorlinks=true,hypertexnames=false, pdfstartview=FitV, linkcolor=blue, citecolor=blue, urlcolor=blue,pagebackref=false]{hyperref}

\usepackage{tikz}
\usepackage{float}
\usepackage[font=footnotesize]{caption}

\usepackage{times,theorem,latexsym,color,comment}

\newcommand{\BOX}{\ensuremath\Box}
\newtheorem{theorema}{Theorem}

\newtheorem{theorem}{Theorem }[section]

{\theorembodyfont{\rmfamily}}
{\theorembodyfont{\rmfamily}}
{\theorembodyfont{\rmfamily}}
\newtheorem{lemma}[theorem]{Lemma}
\newtheorem{proposition}[theorem]{Proposition}
{\theorembodyfont{\rmfamily}\newtheorem{remark}[theorem]{Remark}}
{\theorembodyfont{\rmfamily}}

\newcommand{\Z}{\mathbb{Z}}
\newcommand{\R}{\mathbb{R}}
\newcommand{\C}{\mathbb{C}}
\newcommand{\T}{\mathbb{T}}

\newcommand{\mcL}{\mathcal{L}}

\newcommand{\dd}{\,d}

\newcommand{\eps}{\varepsilon}

\newcommand{\opdist}{\operatorname{dist}}
\newcommand{\ii}{\mathrm i}
\newcommand{\Proj}[1]{\mathcal{P}_{#1}}

\newcommand{\bv}[1]{\bm{#1}} 

\def\XXint#1#2#3{{\setbox0=\hbox{$#1{#2#3}{\int}$}
		\vcenter{\hbox{$#2#3$}}\kern-.5\wd0}}

\newenvironment{proof}{{\vskip\baselineskip\noindent\textbf{Proof:}}}
{\hspace*{.1pt}\hspace*{\fill}\BOX\vskip\baselineskip}

\newenvironment{proofx}[1]
{\vskip\baselineskip\noindent\textbf{Proof of {#1}:}}
{\hspace*{.1pt}\hspace*{\fill}\BOX\vskip\baselineskip}

\begin{document}

    \title{An Arnold-Type Stability Analysis of Vortex Filaments}

\author{
Masashi Aiki
\thanks{
Department of Mathematics, 
Faculty of Science and Technology, 
Tokyo University of Science,
2641 Yamazaki, Noda, Chiba 278-8510, Japan.
\textit{E-mail address:} \texttt{a27120@rs.tus.ac.jp}
}
\and
Mitsuo Higaki
\thanks{
Department of Mathematics, 
Graduate School of Science, 
Kobe University, 
1-1 Rokkodai, Nada-ku, Kobe 657-8501, Japan.
\textit{E-mail address:} \texttt{higaki@math.kobe-u.ac.jp}
}
}
\date{}

\maketitle

\begin{abstract}
We establish nonlinear orbital stability of circular arcs with opening angle $\theta\le \pi$ and closed circles under the Localized Induction Equation (LIE), with prescribed endpoint constraints for arcs. Our Arnold-type variational approach combines two conserved quantities, the bending energy and a component of the fluid impulse, into a functional that we prove is comparable to the squared $H^2$-distance from the corresponding orbit for all admissible curves. For circular arcs with angle $\theta>\pi$, we prove nonlinear orbital instability using the inverse generalized Hasimoto transformation. 
\end{abstract}

\tableofcontents

    \section{Introduction}
    \label{intro}

In this paper, we study the nonlinear stability of circular arcs and closed circles evolving under the following Localized Induction Equation (LIE): 
\begin{equation}
\label{eq.intro.LIE}
\tag{LIE}
    \bv{x}_t=\bv{x}_s\times\bv{x}_{ss}. 
\end{equation}
Here $\bv{x}(s,t)\in\R^3$ parametrizes the evolving curve, $s$ is the arclength parameter, and $t$ denotes time. The LIE is a classical model for the self-induced motion of thin vortex filaments in an inviscid incompressible fluid; see \cite{DaRios1906,Murakamietal1937,ArmsHama1965}. For rigorous results on vortex filaments and rings, see \cite{JerrardSeis2017,GuoJeongZhao2026} for the Euler equations, \cite{GigaMiyakawa1989,BedrossianGermainHarropGriffiths2023,FontelosIspizuaVega2026} for the Navier--Stokes equations, and \cite{delPinoJunemanMussoWei2026} for the Gross--Pitaevskii equation. Circular arcs and closed circles provide basic examples of vortex filaments that translate without changing their shape.

For $R>0$, the LIE admits the following exact uniformly translating solution
\[
    \bv{x}^{R}(s,t)
    :=
    \bv{x}^R_{0}(s)
    + \frac{t}{R}\bv{e}_3,
    \qquad
    \bv{x}^R_{0}(s)
    :=
    R\cos\left(\frac{s}{R}\right)\bv{e}_1
    +
    R\sin\left(\frac{s}{R}\right)\bv{e}_2,
\]
where $\{\bv{e}_1,\bv{e}_2,\bv{e}_3\}$ is the standard basis of $\R^3$. On an interval of length $\theta R$, this solution describes a circular arc of opening angle $\theta$; on the torus of length $2\pi R$, it describes a closed circle. Nonlinear stability estimates for circular arcs and, under symmetry assumptions, closed circles were obtained in \cite{Aiki2025}. Orbital stability of closed circles under non-symmetric perturbations was subsequently established in \cite{AikiHigaki2026a}. Building on these results, we develop a common variational framework for the stable regimes, based on a conserved functional globally comparable to the squared $H^2$-distance from the corresponding orbit. We also identify how the boundary conditions determine the stability and instability.

Our stability arguments are in the spirit of Arnold's variational principle \cite{Arnold1965}: we seek a conserved functional that is coercive modulo the symmetries of the reference filament $\bv{x}^{R}$; see \cite{GallaySverak2024a} for a related application to the stability of planar vortices. To this end, we use the bending energy $B$ (see \cite{LangerPerline1991}) and the third component $P_3$ of the fluid impulse $\bv{P}$, defined by
\[
    B[\bv{x}]
    =
    \frac12\int_D |\bv{x}_{ss}|^2 \dd s,
    \qquad
    P_3[\bv{x}]
    =
    \frac12\int_D
    \bv{e}_3\cdot(\bv{x}\times\bv{x}_s)
    \dd s,
\]
where $D$ denotes the interval or torus under consideration. The fluid impulse and related integral invariants of vortex filaments were studied by Arms and Hama \cite{ArmsHama1965} and Fukumoto \cite{Fukumoto1987}. Both $B$ and $P_3$ are conserved under the LIE evolutions considered here; in the arc problem, the conservation of $P_3$ relies on the endpoint constraints in \eqref{eq.intro.arc.planes} below. We combine $B$ and $P_3$, choosing the coefficient so that $\bv{x}^R_0$ is a critical point under the corresponding boundary-compatible perturbations. This leads to the Lagrange functional
\[
    \mcL[\bv{x}]
    =
    B[\bv{x}] - \frac1{R^3}P_3[\bv{x}],
\]
which also satisfies $\mcL[\bv{x}^R_0]=0$.

The key step in our framework is a quantitative comparison between $\mcL$ and the distance to the orbit of $\bv{x}^{R}$. For circular arcs with
$0<\theta\le\pi$ and for closed circles, we prove
\[
    c\opdist(\bv{x},\mathsf{S})^2
    \le
    \mcL[\bv{x}]
    \le
    C\opdist(\bv{x},\mathsf{S})^2
\]
for all admissible curves $\bv{x}$, where $\opdist(\bv{x},\mathsf{S}):=\inf_{\bv{y}\in\mathsf{S}}\|\bv{x}-\bv{y}\|_{H^2(D)}$, and $c,C>0$ depend only on the fixed geometry. Here $\mathsf{S}$ denotes the orbit $\mathsf{A}_\theta$ or $\mathsf{C}$ defined below, accounting for the allowed translations and, in the circle problem, rotations about the $z$-axis. The admissibility here includes the constraint $|\bv{x}_s|=1$ on the fixed domain $D$, together with the prescribed endpoint constraints in the arc problem or periodicity in the circle problem.

In particular, this comparison requires no smallness assumption on curves $\bv{x}$. Combining it with the conservation of $\mcL[\bv{x}(t)]$ under the LIE evolution immediately yields
\[
    \bigg(
    \frac{c}{C}
    \bigg)^{1/2}
    \opdist(\bv{x}_0,\mathsf{S})
    \le
    \opdist(\bv{x}(t),\mathsf{S})
    \le
    \bigg(
    \frac{C}{c}
    \bigg)^{1/2}
    \opdist(\bv{x}_0,\mathsf{S}),
    \quad 
    t\ge0.
\]
Thus $\mcL$ serves as a Lyapunov functional establishing orbital stability.

By contrast, for circular arcs with $\pi<\theta<2\pi$, the coercivity of $\mcL[\bv{x}]$ fails. This threshold corresponds to the onset of modulational instability of the associated focusing cubic nonlinear Schr\"odinger equation on the interval of length $\theta R$. Using an explicit Akhmediev breather \cite{AkhmedievKorneev1986} and the inverse generalized Hasimoto transformation \cite{Aiki2023}, we can construct solutions of the LIE satisfying the endpoint constraints and prove nonlinear orbital instability. Related constructions for vortex filaments without boundaries were studied in \cite{Salman2013}.

    \subsection{Overview of Our Results}

Since our analysis spans over multiple problem settings, we first give an overview of the main statements for each problem setting. The unique global solutions used in Theorems \ref{thm:intro.arc} and \ref{thm:intro.circ} below are provided by the well-posedness results of \cite{Aiki2017} for the arc problem and \cite{TaniNishiyama1997} for the closed circle problem. Hence, the existence of the solution is not a concern, and only the stability is of importance.

    \paragraph{Circular arcs.}

Let $0<\theta<2\pi$, $L=\theta R$, and $I=(0,L)$. We then consider 
\begin{align}
\label{eq.intro.arc}
\tag{Arc}
    \left\{
    \begin{array}{ll}
    \bv{x}_{t}
    =
    \bv{x}_{s}\times \bv{x}_{ss},
    &(s,t)\in I\times (0,\infty), \\[3mm]
    \bv{x}(s,0)
    =
    \bv{x}_{0}(s),
    &s\in I, \\[3mm]
    \bv{x}_{s}(0,t)
    =
    \bv{e}_{2},
    \quad
    \bv{x}_{s}(L,t)
    =
    \bv{b}_\theta, 
    &t\in (0,\infty), 
    \end{array}\right. 
\end{align}
where $\bv{b}_\theta:=-\sin\theta\,\bv{e}_1+\cos\theta\,\bv{e}_2$. For the initial data, we require that 
\begin{equation}
\label{eq.intro.arc.planes}
    \bv{x}_0(0)\cdot\bv{e}_2
    =
    0,
    \qquad
    \bv{x}_0(L)\cdot\bv{b}_\theta
    =
    0.
\end{equation}
These conditions place the initial endpoints in the same respective planes as those of the arc $\bv{x}^R$. The LIE evolution preserves these constraints, and $\bv{b}_\pi=-\bv{e}_2$. Accordingly, we set
\[
\begin{split}
    \mathsf{A}_{\theta}
    :=
    \begin{cases}
    \bv{x}^R_{0}+\R\bv{e}_3,
    &\theta\ne\pi,\\[3mm]
    \bv{x}^R_{0}+\R\bv{e}_1+\R\bv{e}_3,
    &\theta=\pi,
    \end{cases}
    \qquad
    \opdist(\bv{x},\mathsf{A}_{\theta})
    :=
    \inf_{\bv{y}\in\mathsf{A}_{\theta}}
    \|\bv{x}-\bv{y}\|_{H^2(I)}.
\end{split}
\]

\begin{theorema}
\label{thm:intro.arc}
The following hold.
\begin{enumerate}
\item
\label{item1.thm:intro.arc}
Let $0<\theta\le\pi$. Suppose that $\bv{x}_0\in H^4(I)^3$ satisfies $|\bv{x}_{0,s}|=1$, \eqref{eq.intro.arc.planes}, and the compatibility conditions \eqref{cond1.thm:eq.arc}. Then the corresponding unique global solution of \eqref{eq.intro.arc} satisfies
\[
    c
    \opdist(\bv{x}_0,\mathsf{A}_{\theta})
    \le
    \opdist(\bv{x}(t),\mathsf{A}_{\theta})
    \le 
    C
    \opdist(\bv{x}_0,\mathsf{A}_{\theta}),
    \quad
    t\ge0, 
\]
where the constants $c,C>0$ depend only on $R,\theta$.

\item
\label{item2.thm:intro.arc}
Let $\pi<\theta<2\pi$. Then there exist constants $C_0,\eta>0$, depending only on $R,\theta$, such that, for every $0<\eps<1/2$, there are smooth initial data $\bv{x}^{(\eps)}_{0}$ satisfying $|\bv{x}^{(\eps)}_{0,s}|=1$, \eqref{eq.intro.arc.planes}, and the compatibility conditions \eqref{cond1.thm:eq.arc} with
\[
    \|\bv{x}^{(\eps)}_{0}-\bv{x}^{R}_{0}\|_{H^4(I)}
    \le
    C_0
    \eps
\]
and the corresponding unique global solutions of \eqref{eq.intro.arc} satisfy
\[
    \opdist(\bv{x}^{(\eps)}(T^\eps),\mathsf{A}_\theta)
    \ge
    \eta,
\]
at times $T^\eps$ of order $\log(1/\eps)$.
\end{enumerate}
\end{theorema}

See Section \ref{sec.arc}, especially Theorems \ref{thm:<pi.stab}, \ref{thm:=pi.stab}, and \ref{thm:>pi.instab}, for the precise statements. Theorem \ref{thm:intro.arc} identifies $\theta=\pi$ as the sharp stability threshold. In particular, the semicircle when $\theta=\pi$ is orbitally stable without the symmetry conditions assumed in \cite[Theorem 2.5]{Aiki2025}.

\paragraph{Closed circles.}

Let $L=2\pi R$ and $\T=\R/(L\Z)$. We then consider 
\begin{align}
\label{eq.intro.circ}
\tag{Circ}
    \left\{
    \begin{array}{ll}
    \bv{x}_{t}
    =
    \bv{x}_{s}\times \bv{x}_{ss},
    &(s,t)\in \T\times (0,\infty), \\[3mm]
    \bv{x}(s,0)
    =
    \bv{x}_{0}(s),
    &s\in \T.
    \end{array}\right. 
\end{align}
Let $Q_z^\alpha$ denote the rotation matrix around the $z$-axis by angle $\alpha$. Then we set
\[   
    \mathsf{C}
    := 
    \{
    Q_z^\alpha \bv{x}^R_{0} + \bv{\tau}
    :
    \alpha\in\R,\ \bv{\tau}\in\R^3
    \},
    \qquad
    \opdist(\bv{x},\mathsf{C})
    :=
    \inf_{\bv{y}\in\mathsf{C}}
    \|\bv{x}-\bv{y}\|_{H^2(\T)}.
\]
For the circle $\bv{x}^R_{0}$, rotations about the $z$-axis are equivalent to shifts of the arclength parameter $s$. Thus $\mathsf{C}$ accounts for the phase of the parametrization and spatial translations.

\begin{theorema}
\label{thm:intro.circ}
Suppose that $\bv{x}_0\in H^4(\T)^3$ satisfies $|\bv{x}_{0,s}|=1$. Then the corresponding unique global solution of \eqref{eq.intro.circ} satisfies
\[
    c
    \opdist(\bv{x}_0,\mathsf{C})
    \le
    \opdist(\bv{x}(t),\mathsf{C})
    \le 
    C
    \opdist(\bv{x}_0,\mathsf{C}),
    \quad
    t\ge0,
\]
where the constants $c,C>0$ depend only on $R$.
\end{theorema}

See Section \ref{sec.circ}, especially Theorem \ref{thm:circ.stab}, for the precise statement. Compared with the previous result \cite{AikiHigaki2026a}, Theorem \ref{thm:intro.circ} requires no separate condition on the initial fluid impulse and directly controls the orbital distance $\opdist(\bv{x}(t),\mathsf{C})$ by its initial value.

As a byproduct of the proof of Theorem \ref{thm:intro.circ}, we derive an a priori estimate for the orbital distance of closed filaments under arclength-preserving forcing. This illustrates the usefulness of the same functional even when the conservation laws are no longer available.

    \subsection{Contribution and Related Results}

    \paragraph{Stability of arcs and circles.}

The linear stability of closed vortex filaments was investigated in \cite{CaliniKeithLafortune2011}, where the Hasimoto map is combined with squared eigenfunctions for the nonlinear Schr\"odinger equation. The linear results in \cite{CaliniKeithLafortune2011} are distinct from the nonlinear, uniform-in-time orbital estimates considered in this paper.

For circular arcs, \cite{Aiki2025} established nonlinear stability estimates. When $0<\theta<\pi$, perturbations are controlled uniformly in time in $H^2$ modulo translations along the $z$-axis, without a smallness assumption on the admissible initial perturbation. When $\pi\le\theta<2\pi$, analogous estimates were obtained under an additional midpoint-reflection symmetry.

For closed circles, \cite{AikiHigaki2026a} subsequently established nonlinear orbital stability in $H^2$ under non-symmetric perturbations, modulo spatial translations and rotations about the symmetry axis. Their theorem is formulated in terms of quantitative bounds on the initial fluid impulse relative to that of the reference circle and on the higher Fourier modes. The proof combines a conserved relative energy, which controls the modes $|n|\ge2$, with conservation of the vector fluid impulse, which controls the modes $n=\pm1$.

Building on these results, the present work places the stable arc and circle problems in a common variational framework: the conserved functional $\mcL$ introduced above is globally comparable to the squared $H^2$-distance from the corresponding orbit on the admissible class. In particular, we obtain stability of the semicircle $\theta=\pi$ without reflection symmetry. This approach also allows us to refine the analysis to obtain estimates from below for the orbital distance from the stable arc and circle. Moreover, it allows us to conclude that, if the initial perturbation is not zero, then the solution of the Localized Induction Equation is not asymptotic to the reference arc or circle.

    \paragraph{Variational methods in fluid mechanics.}

Arnold's variational approach \cite{Arnold1965} links the nonlinear stability of steady planar Euler flows to the properties of the kinetic energy under constraints imposed by conserved vorticity quantities. In the energy-Casimir formulation, one supplements the energy with suitable vorticity invariants and establishes coercive estimates for the resulting conserved functional; see \cite{HolmMarsdenRatiuWeinstein1985} for a systematic development of this method in fluid and plasma dynamics. See Arnold \cite{Arnold1966} (also \cite{ArnoldKhesin2021_book}) for a broader account of the geometric and Hamiltonian framework of ideal hydrodynamics.

Variational principles involving the fluid impulse also play an important role in the study of steadily translating Euler vortices. Friedman and Turkington \cite{FriedmanTurkington1981} constructed axisymmetric vortex rings by maximizing kinetic energy under constraints that include the axial impulse. Building on this setting, Choi \cite{Choi2024} proved the orbital stability of Hill's spherical vortex by a concentrated compactness method for axisymmetric flows without swirl. For a recent variational approach to the nonlinear stability of steady
Euler vortex rings, see \cite{CaoLaiQinZhanZou2026}.

For planar vortices, Gallay and {\v{S}}ver\'ak \cite{GallaySverak2024a} investigated the coercivity of quadratic forms associated with Arnold's principle and gave a new stability proof for Oseen's vortex. In \cite{GallaySverak2024b}, they further adapted Arnold-type stability estimates to the vanishing-viscosity analysis of axisymmetric vortex rings. These works demonstrate that variational stability structures can remain useful even when viscosity destroys the inviscid conservation laws. A more recent application to viscous vortex dynamics can also be found in \cite{DolceGallay2026}.

    \paragraph{The Hasimoto transformation and instability.}

Hasimoto \cite{Hasimoto1972} discovered that the evolution of the curvature and torsion under the LIE can be transformed into the focusing cubic nonlinear Schr\"odinger (NLS) equation. Langer and Perline \cite{LangerPerline1991} studied the Poisson geometry of this correspondence and its relation to Hamiltonian flows. Koiso \cite{Koiso1997} developed a moving-frame formulation of that correspondence. For finite-length filaments, Aiki \cite{Aiki2023} established a correspondence between the LIE with endpoint constraints and the focusing cubic NLS equation with homogeneous Neumann boundary conditions. These results provide a framework for studying the stability of vortex filaments through the NLS dynamics.

An explicit nonlinear modulational instability for the focusing NLS equation is provided by the Akhmediev breather \cite{AkhmedievKorneev1986}. Using the correspondence with the LIE, breather-type vortex filaments are constructed in \cite{Salman2013} and their deformations are related to instability. On the NLS side, the instability of the Akhmediev breather is proved in \cite{AlejoFanelliMunoz2019}.

    \paragraph{Quantitative rigidity of almost minimizers.}

Our estimates also admit a static interpretation as quantitative rigidity inequalities. Gallay and {\v{S}}ver\'ak \cite[Remark 2.15]{GallaySverak2024a} emphasize the importance of estimating the distance from almost maximizers to the set of exact maximizers, a question also arising in the stability of the isoperimetric inequality \cite{FuscoMaggiPratelli2008} and of the Sobolev inequality \cite{BianchiEgnell1991}. For recent advances in quantitative stability for sharp functional inequalities, see \cite{DolbeaultEstebanFigalliFrankLoss2025}. In our setting, the corresponding statement concerns minimizers: for circular arcs with $0<\theta\le\pi$ and for closed circles, $\mcL$ has minimum zero on the respective admissible class, and its set of minimizers is $\mathsf{S}$. The lower bound for $\mcL$ yields
\[
    \opdist(\bv{x},\mathsf{S})
    \le
    c^{-1/2}
    \mcL[\bv{x}]^{1/2}.
\]
Thus $\mcL[\bv{x}]^{1/2}$ quantitatively controls the $H^2$-distance to the set $\mathsf{S}$, without a smallness assumption on the admissible curve. A particular feature of our setting is that $\mcL$ is quadratic and its first derivative at $\bv{x}_0^R$ satisfies $D\mcL(\bv{x}_0^R)[\bv{x}-\bv{x}_0^R]=0$ for every admissible $\bv{x}$. Consequently, we have the exact identity $\mcL[\bv{x}]=\mcL[\bv{x}-\bv{x}_0^R]$, with no Taylor remainder. Conservation of $\mcL$ then converts this static quantitative rigidity into nonlinear orbital stability.

\bigskip

The paper is organized as follows. Section \ref{sec.arc} establishes the orbital stability and instability of circular arcs. Section \ref{sec.circ} establishes the orbital stability of closed circles. Section \ref{sec.appl} gives an application to the forced equation.

    \subsection*{Notation}

\begin{itemize}
\item
For $\bv{a},\bv{b}\in\R^3$, the symbols $\bv{a}\cdot\bv{b}$, $\bv{a}\times\bv{b}$, and $|\bv{a}|$ denote the Euclidean inner product, cross product, and norm, respectively. Subscripts $s$ and $t$ indicate partial differentiation. Unless otherwise specified, $c,C$ denote positive constants that may change from line to line; their dependence on parameters is indicated in the relevant statements.

\item
The spatial domain is either an interval $I=(0,L)$ or a torus $\T=\R/(L\Z)$, with $L>0$ specified in each setting. In this paragraph, $D$ stands for either $I$ or $\T$. We use the usual Lebesgue spaces $L^p(D)$, $1\le p\le\infty$, and Sobolev spaces $H^m(D)$, where $m$ is a nonnegative integer. The corresponding spaces of vector-valued functions are denoted by $L^p(D)^3$ and $H^m(D)^3$. When the domain is clear, we abbreviate $\|\bv{u}\|:=\|\bv{u}\|_{L^2(D)}$ and $\|\bv{u}\|_m:=\|\bv{u}\|_{H^m(D)}$, where $\|\bv{u}\|_{H^m(D)}^2:=\sum_{j=0}^{m}\int_D|\partial_s^j\bv{u}(s)|^2\dd s$.

\item
For a $C^2$ functional $F$ on an open subset of a real Hilbert space $H$, we write $F[\bv{x}]$ for its value and denote its first and second Fr\'echet derivatives at $\bv{x}$ by $DF(\bv{x})$ and $D^2F(\bv{x})$, respectively. These are regarded as a linear form and a symmetric bilinear form, whose evaluations on variations are written as $DF(\bv{x})[\bv{\varphi}]$ and $D^2F(\bv{x})[\bv{\varphi},\bv{\psi}]$. For a continuous homogeneous quadratic functional $G:H\to\R$, the second derivative is independent of the evaluation point, and we write simply $D^2G$. In this case,
\[
    G[\bv{x}+\bv{h}]
    =
    G[\bv{x}]
    + DG(\bv{x})[\bv{h}]
    + G[\bv{h}],
    \qquad
    G[\bv{h}]
    =
    \frac12 D^2G[\bv{h},\bv{h}].
\]
\end{itemize}

    \section{Circular Arcs}
    \label{sec.arc}

In this section, we discuss the orbital stability and instability of circular arcs. Let $R>0$ be the radius of the circular arc, $0<\theta<2\pi$ its central angle, $L=\theta R$ its length, and $I=(0,L)$. We then consider the following initial-boundary value problem for the Localized Induction Equation (LIE): 
\begin{align}
\label{eq.arc}
\tag{Arc}
    \left\{
    \begin{array}{ll}
    \bv{x}_{t}
    =
    \bv{x}_{s}\times \bv{x}_{ss},
    &(s,t)\in I\times (0,\infty), \\[3mm]
    \bv{x}(s,0)
    =
    \bv{x}_{0}(s),
    &s\in I, \\[3mm]
    \bv{x}_{s}(0,t)
    =
    \bv{e}_{2},
    \quad
    \bv{x}_{s}(L,t)
    =
    \bv{b}_\theta, 
    &t\in (0,\infty), 
    \end{array}\right. 
\end{align}
where $\bv{x}_{0}$ is the initial filament and
\begin{equation}
\label{def.b}
    \bv{b}_\theta
    :=
    (-\sin\theta,\cos\theta,0).    
\end{equation}
The circular arc $\bv{x}^{R}$ is given by
\[
    \bv{x}^{R}(s,t) 
    := 
    \bv{x}^R_{0}(s)
    + \frac{t}{R}\bv{e}_{3},
    \qquad
    \bv{x}^R_{0}(s)
    := 
    R\cos\left(\frac{s}{R}\right)\bv{e}_1
    + 
    R\sin\left(\frac{s}{R}\right)\bv{e}_2.
\]

We recall the solvability of \eqref{eq.arc}. The following result is a consequence of \cite[Theorem 3.1]{Aiki2023}, applied to the tangent vector $\bv{v}:=\bv{x}_s$, together with the reconstruction of $\bv{x}$.

\begin{theorem}
\label{thm:eq.arc}
Suppose that an initial datum $\bv{x}_0\in H^4(I)^3$ satisfies 
\begin{equation}
\label{cond0.thm:eq.arc}
    |\bv{x}_{0,s}(s)|=1,
    \quad
    s\in\overline{I}
\end{equation}
and the compatibility conditions
\begin{equation}
\label{cond1.thm:eq.arc}
\begin{split}
    &\bv{x}_{0,s}(0)
    =
    \bv{e}_2,
    \qquad
    \bv{x}_{0,s}(L)
    =
    \bv{b}_\theta, \\[3mm]
    &(\bv{x}_{0,s}\times\bv{x}_{0,sss})(0)
    =
    (\bv{x}_{0,s}\times\bv{x}_{0,sss})(L)
    =
    \bv{0}.
\end{split}
\end{equation}
Then \eqref{eq.arc} has a unique global-in-time solution $\bv{x}$ in the class
\begin{equation}
\label{est1.thm:eq.arc}
    \bv{x}
    \in
    \bigcap_{j=0}^{2}
    C^j([0,\infty);H^{4-2j}(I)^3).
\end{equation}
Moreover, $\bv{x}$ satisfies 
\begin{equation}
\label{est2.thm:eq.arc}
\begin{split}
    &\bv{x}_s
    \in 
    L^\infty(0,\infty; H^3(I)^3), \\[3mm]
    &|\bv{x}_s(s,t)|=1,
    \quad
    (s,t)\in\overline{I}\times[0,\infty),    
\end{split}
\end{equation}
and for all $t\in [0,\infty)$, 
\begin{equation}
\label{est3.thm:eq.arc}
\begin{split}
    \bv{x}(0,t)\cdot\bv{e}_2
    =
    \bv{x}_0(0)\cdot\bv{e}_2, 
    \qquad
    \bv{x}(L,t)\cdot\bv{b}_\theta
    =
    \bv{x}_0(L)\cdot\bv{b}_\theta.
\end{split}
\end{equation}
\end{theorem}

    \subsection{Variational Setting for Circular Arcs}
    \label{sec.arc.setting}

We introduce the variational framework for the arc problem \eqref{eq.arc}. The following definitions and lemmas are valid for all $0<\theta <2\pi$.

    \subsubsection*{Function Spaces}

Define 
\begin{equation}
    X_{\theta}
    :=
    \bigg\{
    \bv{\varphi}\in H^2(I)^3 
    :
    \begin{array}{l}
    \bv{\varphi}_s(0)
    =
    \bv{\varphi}_s(L)
    =
    \bv{0}, \\[3pt]
    \bv{\varphi}(0)\cdot\bv{e}_2
    =
    \bv{\varphi}(L)\cdot\bv{b}_\theta
    =
    0
    \end{array}
    \bigg\}
\end{equation}
and
\begin{equation}
\label{eq:mcM}
    Y_{\theta}
    :=
    \bv{x}^R_{0}
    + X_{\theta}
    = 
    \{
    \bv{x}^R_{0} + \bv{\varphi}
    :
    \bv{\varphi} \in X_{\theta}
    \}.
\end{equation}
Finally, we set
\begin{equation}
    M_{\theta}
    = 
    \{
    \bv{x}\in Y_\theta 
    :
    |\bv{x}_s(s)|=1\ \text{for all }s\in I
    \}.
\end{equation}

    \subsubsection*{Lagrange Functional}

For $\bv{x}\in H^2(I)^3$, we define
\[
    B[\bv{x}]
    := 
    \frac12
    \int_{I}
    |\bv{x}_{ss}|^2
    \dd s
\]
and
\[
    P_3[\bv{x}]
    :=
    \frac12
    \int_{I}
    \bv{e}_3
    \cdot
    (\bv{x}\times\bv{x}_s)
    \dd s.
\]
A direct computation shows that 
\begin{equation}
\label{rep.P3}
    P_3[\bv{x}]
    =
    \frac12
    \int_I
    (x_1 x_{2,s}
    - x_2 x_{1,s})
    \dd s.
\end{equation}

Consider the following Lagrange functional $\mcL$ on $H^2(I)^3$: 
\begin{equation}
\label{def.L.arc}
\begin{split}
    \mcL[\bv{x}]
    &:=
    B[\bv{x}]
    -\frac1{R^3}
    P_3[\bv{x}] \\[3mm]
    &=
    \frac12
    \bigg(
    \int_{I}|\bv{x}_{ss}|^2\dd s
    -\frac1{R^3}\int_{I}
    \bv{e}_3\cdot(\bv{x}\times\bv{x}_s)\dd s
    \bigg).
\end{split}    
\end{equation}
The functional $\mcL$ is smooth and quadratic on $H^2(I)^3$, and its derivative is
\begin{equation}
\label{eq.DL.arc}
    D\mcL(\bv{x})[\bv{\varphi}]
    =
    \int_I
    \bv{x}_{ss}\cdot\bv{\varphi}_{ss}
    \dd s 
    - \frac1{2R^3}
    \int_I
    \bv{e}_3\cdot
    (
    \bv{\varphi}\times\bv{x}_s
    + \bv{x}\times\bv{\varphi}_s
    )
    \dd s.
\end{equation}
Note that $\mcL$ is invariant under translations along $\bv{e}_3$: for any $\bv{x}\in H^2(I)^3$ and $a\in\R$, 
\begin{equation}
\label{eq.L.arc.inv}
    \mcL[\bv{x} + a\bv{e}_3]
    =
    \mcL[\bv{x}].
\end{equation}
Moreover, for any $0<\theta<2\pi$, 
\begin{equation}
\label{eq.L.arc.xR0}
    B[\bv{x}^R_{0}]
    =
    \frac{\theta}{2R},
    \qquad
    P_3[\bv{x}^R_{0}]
    =
    \frac{\theta R^2}{2},
    \qquad
    \mcL[\bv{x}^R_{0}]
    =
    0.
\end{equation}

The following lemma shows that $\mcL[\bv{x}]$ is conserved for the solutions in Theorem \ref{thm:eq.arc}.

\begin{lemma}
\label{lem:L.arc.LOC}
Let $\bv{x}=\bv{x}(t)$ be the unique global solution of \eqref{eq.arc} in Theorem \ref{thm:eq.arc} with the initial datum $\bv{x}_0\in M_\theta\cap H^4(I)^3$ satisfying \eqref{cond1.thm:eq.arc}. Then $\bv{x}(t)\in M_{\theta}\cap H^4(I)^3$ holds for all $t\ge0$. Moreover, $B[\bv{x}(t)]$, $P_3[\bv{x}(t)]$, and $\mcL[\bv{x}(t)]$ are conserved in time.
\end{lemma}

\begin{proof}
The assertion $\bv{x}(t)\in M_{\theta}$ for $t\ge0$ follows from Theorem \ref{thm:eq.arc}. The regularity in \eqref{est1.thm:eq.arc} justifies the time differentiations and integrations by parts below. Differentiating the equation in \eqref{eq.arc} with respect to $s$, and the boundary conditions with respect to $t$, we have
\[
    \bv{x}_{st}
    =
    \bv{x}_s\times\bv{x}_{sss},
    \qquad
    \bv{x}_{st}(0,t)
    =
    \bv{x}_{st}(L,t)
    =
    \bv{0}.
\]
Then the conservation of $B[\bv{x}(t)]$ follows from
\begin{align*}
    \frac{d}{dt}
    B[\bv{x}(t)]
    &=
    \int_I 
    \bv{x}_{ss}\cdot\bv{x}_{sst}
    \dd s \\
    &=
    \big[
    \bv{x}_{ss}\cdot\bv{x}_{st}
    \big]_{s=0}^{s=L}
    -\int_I
    \bv{x}_{sss}\cdot(\bv{x}_s\times\bv{x}_{sss})
    \dd s
    =0.
\end{align*}

We next prove the conservation of $P_3[\bv{x}(t)]$. Since $\bv{x}_0\in Y_\theta$, we have $\bv{x}_0(0)\cdot\bv{e}_2=\bv{x}_0(L)\cdot\bv{b}_\theta=0$. Thus, by \eqref{est3.thm:eq.arc} and the boundary conditions, at both $s=0,L$, 
\[
    \bv{x}\cdot\bv{x}_s
    =
    \bv{e}_3\cdot\bv{x}_s
    =
    0.
\]
It follows that, at both $s=0,L$,
\begin{align*}
    \bv{e}_3\cdot(\bv{x}\times\bv{x}_t)
    &=
    \bv{e}_3\cdot
    (\bv{x}\times(\bv{x}_s\times\bv{x}_{ss})) \\[3mm]
    &=
    (\bv{e}_3\cdot\bv{x}_s)(\bv{x}\cdot\bv{x}_{ss})
    -
    (\bv{e}_3\cdot\bv{x}_{ss})(\bv{x}\cdot\bv{x}_s)
    =0.
\end{align*}
Moreover, $|\bv{x}_s|=1$ gives $\bv{x}_s\cdot\bv{x}_{ss}=0$, and hence
\[
    \bv{x}_t\times\bv{x}_s
    =
    (\bv{x}_s\times\bv{x}_{ss})\times\bv{x}_s
    =
    \bv{x}_{ss}.
\]
Therefore, integration by parts yields
\begin{align*}
    \frac{d}{dt}
    P_3[\bv{x}(t)]
    &=
    \frac12
    \int_I
    \bv{e}_3\cdot
    (\bv{x}_t\times\bv{x}_s+\bv{x}\times\bv{x}_{st})
    \dd s \\[3mm]
    &=
    \int_I
    \bv{e}_3\cdot(\bv{x}_t\times\bv{x}_s)
    \dd s
    + \frac12
    \big[
    \bv{e}_3\cdot(\bv{x}\times\bv{x}_t)
    \big]_{s=0}^{s=L} \\[3mm]
    &=
    \bv{e}_3\cdot
    (\bv{x}_s(L,t)-\bv{x}_s(0,t))
    =0.
\end{align*}
Since $\mcL=B-R^{-3}P_3$, the conservation of $\mcL[\bv{x}(t)]$ follows as well.
\end{proof}

Next we show that the circular arc $\bv{x}^R_{0}$ is a critical point of $\mcL$ restricted to $Y_{\theta}$.

\begin{lemma}
\label{lem:L.arc.deriv}
For any $\bv{\varphi}\in X_{\theta}$, we have
\[
    D\mcL(\bv{x}^R_{0})[\bv{\varphi}]
    =
    0.
\]
\end{lemma}

\begin{proof}
Let $\bv{\varphi}\in X_\theta$. The definition of $X_\theta$ gives 
\[
    \bv{\varphi}_s(0)
    =
    \bv{\varphi}_s(L)
    =
    \bv{0},
    \qquad
    (\bv{x}^R_{0,s}\cdot\bv{\varphi})(0)
    =
    (\bv{x}^R_{0,s}\cdot\bv{\varphi})(L)
    =
    0.
\]
Using
$\bv{x}^R_{0,sss}=-R^{-2}\bv{x}^R_{0,s}$ and
$\bv{x}^R_{0,ssss}=R^{-4}\bv{x}^R_{0}$,
we integrate by parts twice to obtain
\begin{align*}
    \int_I
    \bv{x}^R_{0,ss}\cdot\bv{\varphi}_{ss}
    \dd s
    &=
    \big[
    \bv{x}^R_{0,ss}\cdot\bv{\varphi}_s
    - \bv{x}^R_{0,sss}\cdot\bv{\varphi}
    \big]_{s=0}^{s=L}
    + \int_I
    \bv{x}^R_{0,ssss}\cdot\bv{\varphi}\dd s \\[3mm]
    &=
    \frac1{R^4}
    \int_I
    \bv{x}^R_{0}\cdot\bv{\varphi}
    \dd s.
\end{align*}
Also, since
\[
    \bv{x}^R_{0,s}\times\bv{e}_3
    =\frac1R\bv{x}^R_{0},
    \qquad
    \bv{e}_3\times\bv{x}^R_{0}
    =R\bv{x}^R_{0,s},
\]
we have
\begin{align*}
    &\frac12\int_I
    \bv{e}_3\cdot
    (
    \bv{\varphi}\times\bv{x}^R_{0,s}
    + \bv{x}^R_{0}\times\bv{\varphi}_s
    )
    \dd s \\[3mm]
    &=
    \int_I
    \bv{e}_3\cdot
    (\bv{\varphi}\times\bv{x}^R_{0,s})
    \dd s
    + \frac12
    \big[\bv{e}_3\cdot
    (\bv{x}^R_{0}\times\bv{\varphi})
    \big]_{s=0}^{s=L} \\[3mm]
    &=
    \frac1R
    \int_I
    \bv{x}^R_{0}\cdot\bv{\varphi}
    \dd s.
\end{align*}
Substituting these identities into the formula \eqref{eq.DL.arc} for $D\mcL$ gives
\[
    D\mcL(\bv{x}^R_{0})[\bv{\varphi}]
    =
    \frac1{R^4}\int_I\bv{x}^R_{0}\cdot\bv{\varphi}\dd s
    -\frac1{R^4}\int_I\bv{x}^R_{0}\cdot\bv{\varphi}\dd s
    =0.
\]
This completes the proof.
\end{proof}

Lemma \ref{lem:L.arc.deriv} allows us to express $\mcL[\bv{x}]$ in terms of the perturbation from $\bv{x}^R_{0}$.

\begin{lemma}
\label{lem:L.arc.reduct}
For any $\bv{x}=\bv{x}^R_{0} + \bv{\varphi}\in Y_\theta$ where $\bv{\varphi}\in X_\theta$, we have
\[
    \mcL[\bv{x}] 
    = 
    \mcL[\bv{\varphi}].
\]
\end{lemma}

\begin{proof}
Since $\mcL$ is quadratic, we have 
\[
    \mcL[\bv{x}]
    =
    \mcL[\bv{x}^R_{0}]
    + D\mcL(\bv{x}^R_{0})[\bv{\varphi}]
    + \mcL[\bv{\varphi}].
\]
Then the assertion follows from $\mcL[\bv{x}^R_{0}]=0$ and $D\mcL(\bv{x}^R_{0})[\bv{\varphi}]=0$ by Lemma \ref{lem:L.arc.deriv}. 
\end{proof}

From here on, we divide the analysis according to the size of the central angle $\theta$.

    \subsection{Stable Arc: The Case \texorpdfstring{$\theta<\pi$}{theta < pi}}
    \label{sec.<pi.stab}

Let $0<\theta<\pi$. Define
\[
    \mathsf{A}_{\theta}
    :=
    \bv{x}^R_{0} + \R \bv{e}_3,
    \qquad
    \opdist(\bv{x},\mathsf{A}_{\theta})
    :=
    \inf_{\bv{y}\in\mathsf{A}_{\theta}}
    \|\bv{x}-\bv{y}\|_2.
\]
The aim of this section is to prove the following stability result.

\begin{theorem}
\label{thm:<pi.stab}
Let $\bv{x}=\bv{x}(t)$ be the unique global solution of \eqref{eq.arc} in Theorem \ref{thm:eq.arc} with the initial datum $\bv{x}_0\in M_\theta\cap H^4(I)^3$ satisfying \eqref{cond1.thm:eq.arc}. Then $\bv{x}(t)$ satisfies 
\[
    \bigg(
    \frac{c}{C}
    \bigg)^{1/2}
    \opdist(\bv{x}_0,\mathsf{A}_{\theta})
    \le
    \opdist(\bv{x}(t),\mathsf{A}_{\theta})
    \le
    \bigg(
    \frac{C}{c}
    \bigg)^{1/2}
    \opdist(\bv{x}_0,\mathsf{A}_{\theta}),
    \quad
    t\ge0,
\]
where $c,C>0$ are the constants in Proposition \ref{prop:<pi.equiv} below.
\end{theorem}

We prove Theorem \ref{thm:<pi.stab} as a corollary of the following proposition.

\begin{proposition}
\label{prop:<pi.equiv}
There exist constants $c,C>0$ depending on $R,\theta$ such that, for any $\bv{x}\in Y_\theta$, 
\[
    c\opdist(\bv{x},\mathsf{A}_{\theta})^2
    \le
    \mcL[\bv{x}]
    \le
    C\opdist(\bv{x},\mathsf{A}_{\theta})^2.
\]
In particular, 
\begin{equation}
\label{eq.L.<pi.zero}
    \{
    \bv{x}\in Y_\theta
    :
    \mcL[\bv{x}]=0
    \}
    =
    \mathsf{A}_{\theta}.
\end{equation}
\end{proposition}

To prove Proposition \ref{prop:<pi.equiv}, we need several auxiliary lemmas.

\begin{lemma}
\label{lem:<pi.coercive}
Let $0<\theta<\pi$. For any $\bv{\varphi}\in X_\theta$, we have
\begin{equation}
\label{est.lem:<pi.coercive}
\begin{split}
    \mcL[\bv{\varphi}]
    \ge
    \frac12
    \left(
    1-\frac{\theta^2}{\pi^2}
    \right)
    \|\bv{\varphi}_{ss}\|^2.
\end{split}
\end{equation}
\end{lemma}

\begin{proof}
We first estimate the impulse term $P_3[\bv{\varphi}]$ in \eqref{def.L.arc}. Since
\[
    \varphi_2(0)
    =
    0,
    \qquad
    \varphi_1(L)
    =
    \frac{\cos\theta}{\sin\theta}
    \varphi_2(L), 
\]
integration by parts in \eqref{rep.P3} gives
\[
    P_3[\bv{\varphi}]
    =
    \frac12
    \frac{\cos\theta}{\sin\theta}
    \varphi_2(L)^2 
    - \int_I
    \varphi_{1,s}\varphi_2
    \dd s.
\]
We also have
\begin{equation}
\label{est1.prf.lem:<pi.coercive}
    \frac{\cos\theta}{\sin\theta}\varphi_2(L)^2
    \le
    R\|\varphi_{2,s}\|^2
    -\frac1R\|\varphi_2\|^2.
\end{equation}
Indeed, integration by parts, together with $L=\theta R$, gives
\[
\begin{aligned}
    &R\|\varphi_{2,s}\|^2
    -\frac1R\|\varphi_2\|^2
    -\frac{\cos\theta}{\sin\theta}\varphi_2(L)^2 \\[3mm]
    &=
    R\int_I
    \left(
    \varphi_{2,s}(s)
    -\frac1R
    \frac{\cos(s/R)}{\sin(s/R)}
    \varphi_2(s)
    \right)^2
    \dd s
    \ge 0.
\end{aligned}
\]
The boundary term at $s=0$ vanishes since
$\varphi_2\in H^2(I)$ and $\varphi_2(0)=0$ imply
$\varphi_2(s)=O(s)$ as $s\to0$. 
Then the Cauchy--Schwarz and Young inequalities give 
\[
\begin{split}
    P_3[\bv{\varphi}]
    &\le
    \frac12
    \frac{\cos\theta}{\sin\theta}
    \varphi_2(L)^2
    +
    \|\varphi_{1,s}\|\,\|\varphi_2\| \\[3mm]
    &\le
    \frac12
    \frac{\cos\theta}{\sin\theta}
    \varphi_2(L)^2
    +
    \frac R2\|\varphi_{1,s}\|^2
    +
    \frac1{2R}\|\varphi_2\|^2 \\[3mm]
    &\le
    \frac R2
    (
    \|\varphi_{1,s}\|^2+\|\varphi_{2,s}\|^2
    ).
\end{split}
\]
Consequently,
\begin{equation}
\label{est2.prf.lem:<pi.coercive}
    P_3[\bv{\varphi}]
    \le
    \frac R2\|\bv{\varphi}_s\|^2.    
\end{equation}

The Poincar\'e inequality yields
\[
    \|\bv{\varphi}_s\| 
    \le 
    \frac{L}{\pi}
    \|\bv{\varphi}_{ss}\|.
\]
Hence, using \eqref{est2.prf.lem:<pi.coercive}, we see that \eqref{est.lem:<pi.coercive} follows from 
\[
    \mcL[\bv{\varphi}]
    \ge
    \frac12
    \left(
    \|\bv{\varphi}_{ss}\|^2
    -\frac1{R^2}\|\bv{\varphi}_s\|^2
    \right)
    \ge
    \frac12
    \left(
    1-\frac{\theta^2}{\pi^2}
    \right)
    \|\bv{\varphi}_{ss}\|^2.
\]
The proof is complete.
\end{proof}

The following lemma links the orbital distance and the functional $\mcL$ via the $H^2$-norm of the perturbation.

\begin{lemma}
\label{lem:<pi.proj}
Let $0<\theta<\pi$. For $\bv{x}\in Y_\theta$, let $\bv{\varphi}$ be defined by
\begin{equation}
\label{eq1.lem:<pi.proj}
    \bv{\varphi}
    :=
    \bv{x} - \bv{x}^R_{0} 
    - p \bv{e}_3,
    \qquad
    p
    :=
    \frac{1}{L}
    \int_{I} 
    \bv{e}_3
    \cdot
    (\bv{x} - \bv{x}^R_{0}). 
\end{equation}
Then $\bv{\varphi}$ satisfies the following: 
\begin{enumerate}[(i)]
\item
\label{item1.lem:<pi.proj}
$\bv{\varphi}\in X_\theta$ and $\int_I \bv{e}_3\cdot\bv{\varphi}=0$.

\item
\label{item2.lem:<pi.proj}
$c\|\bv{\varphi}\|_2 \le \|\bv{\varphi}_{ss}\|\le \|\bv{\varphi}\|_2$ with a constant $c$ depending on $R,\theta$.

\item
\label{item3.lem:<pi.proj}
$\opdist(\bv{x},\mathsf{A}_{\theta}) = \|\bv{\varphi}\|_2$.

\item
\label{item4.lem:<pi.proj}
$\mcL[\bv{x}] = \mcL[\bv{\varphi}]$.
\end{enumerate}
\end{lemma}

\begin{proof}
(\ref{item1.lem:<pi.proj})
Using the property of $\bv{x}^R_{0,s}$ and $\bv{e}_3\cdot\bv{e}_2=\bv{e}_3\cdot\bv{b}_\theta=0$, we see that
\[
    \bv{\varphi}_s(0)
    =
    \bv{\varphi}_s(L)
    =
    \bv{0},
    \qquad
    \bv{\varphi}(0)\cdot\bv{e}_2
    =
    \bv{\varphi}(L)\cdot\bv{b}_\theta
    =
    0.
\]
Moreover, $\int_I \bv{e}_3\cdot\bv{\varphi}=0$ by the definition of $p$. This proves (\ref{item1.lem:<pi.proj}).

(\ref{item2.lem:<pi.proj})
The upper bound is obvious. Since $\sin\theta\ne0$, the conditions at $s=0,L$ imply
\[
    \varphi_1(0)
    =
    \frac{\cos\theta}{\sin\theta}
    \int_I\varphi_{2,s}\dd s
    -
    \int_I\varphi_{1,s}\dd s,
    \qquad
    \varphi_2(0)
    =
    0.
\]
Then the fundamental theorem of calculus and the Cauchy--Schwarz inequality yield
\[
    \|(\varphi_1,\varphi_2)\|
    \le
    C \|(\varphi_{1,s},\varphi_{2,s})\|.
\]
In addition, since $\varphi_3$ has zero mean and $\bv{\varphi}_{s}\in H^1_0(I)^3$, 
\[
    \|\varphi_3\| 
    \le 
    \frac{L}{\pi}
    \|\varphi_{3,s}\|,
    \qquad
    \|\bv{\varphi}_s\|
    \le
    \frac{L}{\pi}
    \|\bv{\varphi}_{ss}\|.
\]
Combining these estimates gives the lower bound. This proves (\ref{item2.lem:<pi.proj}).

(\ref{item3.lem:<pi.proj})
For every $a\in\R$, the zero mean of $\bv{e}_3\cdot\bv{\varphi}$ implies
\[
    \|\bv{x}-(\bv{x}^R_{0}+a\bv{e}_3)\|_2^2
    =
    \|\bv{\varphi}+(p-a)\bv{e}_3\|_2^2
    =
    \|\bv{\varphi}\|_2^2
    + L(p-a)^2.
\]
Minimization in $a$ proves (\ref{item3.lem:<pi.proj}).

(\ref{item4.lem:<pi.proj})
This follows from \eqref{eq.L.arc.inv}, (\ref{item1.lem:<pi.proj}) above, and Lemma \ref{lem:L.arc.reduct}. 
\end{proof}

Now we prove Proposition \ref{prop:<pi.equiv} as a corollary of Lemmas \ref{lem:<pi.coercive} and \ref{lem:<pi.proj}.

\begin{proofx}{Proposition \ref{prop:<pi.equiv}}
Let $\bv{x}\in Y_\theta$ and define $\bv{\varphi}\in X_\theta$ by \eqref{eq1.lem:<pi.proj}. By Lemma \ref{lem:<pi.proj} (\ref{item4.lem:<pi.proj}), we have $\mcL[\bv{x}]=\mcL[\bv{\varphi}]$. In view of Lemma \ref{lem:<pi.proj} (\ref{item2.lem:<pi.proj}) and (\ref{item3.lem:<pi.proj}), it suffices to show that 
\[
    c\|\bv{\varphi}_{ss}\|^2
    \le
    \mcL[\bv{x}]
    \le
    C\|\bv{\varphi}\|_2^2, 
\]
for some constants $c,C$ depending on $R,\theta$. The upper bound follows from the definition of $\mcL$ in \eqref{def.L.arc}. The lower bound follows from Lemma \ref{lem:<pi.coercive}. The proof is complete.
\end{proofx}

Combining Proposition \ref{prop:<pi.equiv} with Lemma \ref{lem:L.arc.LOC}, we prove  Theorem \ref{thm:<pi.stab}.

\begin{proofx}{Theorem \ref{thm:<pi.stab}}
Thanks to Lemma \ref{lem:L.arc.LOC}, the solution $\bv{x}=\bv{x}(t)$ satisfies $\bv{x}(t)\in M_\theta\subset Y_\theta$ and $\mcL[\bv{x}(t)]=\mcL[\bv{x}_0]$ for all $t\ge0$. Thus, applying Proposition \ref{prop:<pi.equiv}, we obtain
\[
    c\opdist(\bv{x}(t),\mathsf{A}_{\theta})^2
    \le
    \mcL[\bv{x}(t)] 
    =
    \mcL[\bv{x}_0] 
    \le
    C\opdist(\bv{x}_0,\mathsf{A}_{\theta})^2,
    \quad
    t\ge0.
\]
Dividing by $c$ and taking square roots gives the upper bound. The lower bound is obtained in a similar manner. This completes the proof.
\end{proofx}

    \subsection{Stable Arc: The Case \texorpdfstring{$\theta=\pi$}{theta = pi}}
    \label{sec.=pi.stab}

Compared to the case $\theta < \pi$, the extra translation invariance along the $\bv{e}_{1}$-axis motivates us to define the orbit and orbital distance by the following.
\[
    \mathsf{A}_{\pi}
    :=
    \bv{x}^R_{0} + \R \bv{e}_1 + \R \bv{e}_3,
    \qquad
    \opdist(\bv{x},\mathsf{A}_{\pi})
    := 
    \inf_{\bv{y}\in\mathsf{A}_{\pi}}
    \|\bv{x}-\bv{y}\|_2.
\]
The aim of this section is to prove the following stability result.

\begin{theorem}
\label{thm:=pi.stab}
Let $\bv{x}=\bv{x}(t)$ be the unique global solution of \eqref{eq.arc} in Theorem \ref{thm:eq.arc} with the initial datum $\bv{x}_0\in M_\pi\cap H^4(I)^3$ satisfying \eqref{cond1.thm:eq.arc}. Then $\bv{x}(t)$ satisfies 
\[
    \bigg(
    \frac{c}{C}
    \bigg)^{1/2}
    \opdist(\bv{x}_0,\mathsf{A}_{\pi})
    \le
    \opdist(\bv{x}(t),\mathsf{A}_{\pi})
    \le
    \bigg(
    \frac{C}{c}
    \bigg)^{1/2}
    \opdist(\bv{x}_0,\mathsf{A}_{\pi}),
    \quad
    t\ge0,
\]
where $c,C>0$ are the constants in Proposition \ref{prop:=pi.equiv} below.
\end{theorem}

The proof of Theorem \ref{thm:=pi.stab} is omitted, as it follows an argument similar to that of Theorem \ref{thm:<pi.stab} by combining Lemma \ref{lem:L.arc.LOC} with the following proposition.

\begin{proposition}
\label{prop:=pi.equiv}
There exist constants $c,C>0$  depending on $R$ such that, for any $\bv{x}\in Y_\pi$, 
\[
    c\opdist(\bv{x}, \mathsf{A}_{\pi})^2
    \le
    \mcL[\bv{x}]
    \le
    C\opdist(\bv{x}, \mathsf{A}_{\pi})^2.
\]
In particular, 
\begin{equation}
\label{eq.L.=pi.zero}
    \{
    \bv{x}\in Y_\pi
    :
    \mcL[\bv{x}]=0
    \}
    =
    \mathsf{A}_{\pi}.
\end{equation}
\end{proposition}

\begin{remark}
\label{rem:prop:=pi.equiv}
The proof below gives the choices
\[
    c
    =
    \frac1{4(1+R^2+R^4)}, 
    \qquad
    C
    =
    \frac34.
\]
These constants are not asserted to be optimal.
\end{remark}

To prove Proposition \ref{prop:=pi.equiv}, we need several auxiliary lemmas.

\begin{lemma}
\label{lem:=pi.inv}
$\mcL$ is invariant under translations along $\bv{e}_1$: for any $\bv{x}\in Y_\pi$ and $b\in\R$,
\begin{equation}
\label{eq.L.arc.=pi.inv}
    \mcL[\bv{x} + b\bv{e}_1]
    =
    \mcL[\bv{x}].
\end{equation}
\end{lemma}

\begin{proof}
Since $\bv{x}\in Y_\pi$, we have $x_2(0)=x_2(L)=0$. Thus
\[
    P_3[\bv{x} + b\bv{e}_1]
    -
    P_3[\bv{x}]
    =
    \frac{b}{2}
    (
    x_2(L) - x_2(0)
    )
    =
    0.
\]
Hence both $B$ and $P_3$ are invariant under translations along $\bv{e}_1$, and so is $\mcL$.
\end{proof}

\begin{lemma}
\label{lem:=pi.coercive}
For any $\bv{\varphi}\in X_\pi$, we have
\begin{equation}
\label{est.lem:=pi.coercive}
\begin{split}
    \mcL[\bv{\varphi}]
    \ge
    \frac14
    \|\bv{\varphi}_{ss}\|^2.
\end{split}
\end{equation}
\end{lemma}

\begin{proof}
Since $\varphi_2(0)=\varphi_2(L)=0$,  integration by parts in \eqref{rep.P3} gives
\[
    P_3[\bv{\varphi}]
    =
    \frac12
    \int_I
    (\varphi_1\varphi_{2,s}
    -\varphi_2\varphi_{1,s})
    \dd s
    =
    -\int_I
    \varphi_{1,s}\varphi_2
    \dd s.
\]
It follows that
\begin{equation}
\label{eq1.prf.lem:=pi.coercive}
    \mcL[\bv{\varphi}]
    =
    \frac12\|\bv{\varphi}_{ss}\|^2
    + \frac1{R^3}\int_I\varphi_{1,s}\varphi_2\dd s.
\end{equation}

The Poincar\'e inequality yields
\[
    \|\varphi_{1,s}\|
    \le 
    R\|\varphi_{1,ss}\|,
    \qquad
    \|\varphi_2\|
    \le 
    R\|\varphi_{2,s}\|.
\]
In addition, we have $\varphi_{2,s}(0)=\varphi_{2,s}(L)=0$ and
\[
    \int_I\varphi_{2,s}\dd s
    =
    \varphi_2(L)-\varphi_2(0)
    =
    0.
\]
Thus the periodic extension of $\varphi_{2,s}$ belongs to $H^1$ with zero mean. Hence, 
\[
    \|\varphi_{2,s}\|
    \le
    \frac{L}{2\pi}
    \|\varphi_{2,ss}\|
    =
    \frac{R}{2}
    \|\varphi_{2,ss}\|,
    \qquad
    \|\varphi_2\|
    \le
    \frac{R^2}{2}
    \|\varphi_{2,ss}\|.
\]
Consequently,
\[
\begin{split}
    \frac1{R^3}
    \left|
    \int_I
    \varphi_{1,s}\varphi_2\dd s
    \right|
    \le
    \frac12
    \|\varphi_{1,ss}\|\,
    \|\varphi_{2,ss}\| 
    \le
    \frac14
    (
    \|\varphi_{1,ss}\|^2+\|\varphi_{2,ss}\|^2
    ),
\end{split}
\]
which combined with \eqref{eq1.prf.lem:=pi.coercive} leads to the desired estimate \eqref{est.lem:=pi.coercive}. 
\end{proof}

\begin{lemma}
\label{lem:=pi.proj}
For $\bv{x}\in Y_\pi$, let $\bv{\varphi}$ be defined by
\begin{equation}
\label{eq1.lem:=pi.proj}
    \bv{\varphi}
    :=
    \bv{x} - \bv{x}^R_{0} 
    - p \bv{e}_3
    - q \bv{e}_1, 
\end{equation}
where
\begin{equation}
\label{eq2.lem:=pi.proj}
    p
    :=
    \frac{1}{L}
    \int_{I} 
    \bv{e}_3
    \cdot
    (\bv{x} - \bv{x}^R_{0}),
    \qquad
    q
    :=
    \frac{1}{L}
    \int_{I} 
    \bv{e}_1
    \cdot
    (\bv{x} - \bv{x}^R_{0}). 
\end{equation}
Then $\bv{\varphi}$ satisfies the following: 
\begin{enumerate}[(i)]
\item
\label{item1.lem:=pi.proj}
$\bv{\varphi}\in X_\pi$, $\int_I \bv{e}_3\cdot\bv{\varphi}=0$, and $\int_I \bv{e}_1\cdot\bv{\varphi}=0$.

\item
\label{item2.lem:=pi.proj}
$c\|\bv{\varphi}\|_2 \le \|\bv{\varphi}_{ss}\|\le \|\bv{\varphi}\|_2$ with the constant $c=(1+R^2+R^4)^{-1/2}$.

\item
\label{item3.lem:=pi.proj}
$\opdist(\bv{x},\mathsf{A}_{\pi}) = \|\bv{\varphi}\|_2$.

\item
\label{item4.lem:=pi.proj}
$\mcL[\bv{x}] = \mcL[\bv{\varphi}]$.
\end{enumerate}
\end{lemma}

\begin{proof}
(\ref{item1.lem:=pi.proj})
Using the property of $\bv{x}^R_{0,s}$ and $\bv{e}_j\cdot\bv{e}_2=\bv{e}_j\cdot\bv{b}_\pi=0$
for $j=1,3$, we see that
\[
    \bv{\varphi}_s(0)
    =
    \bv{\varphi}_s(L)
    =
    \bv{0},
    \qquad
    \bv{\varphi}(0)\cdot\bv{e}_2
    =
    \bv{\varphi}(L)\cdot\bv{b}_\pi
    =
    0.
\]
Moreover, $\int_I\bv{e}_3\cdot\bv{\varphi}=\int_I \bv{e}_1\cdot\bv{\varphi}=0$ by the definitions of $p,q$. This proves (\ref{item1.lem:=pi.proj}).

(\ref{item2.lem:=pi.proj})
The upper bound is obvious. Since $\bv{b}_\pi=-\bv{e}_2$, the conditions at $s=0,L$ imply
\[
    \varphi_2(0)
    =
    \varphi_2(L)
    =
    0.
\]
Then, since $\varphi_1$ and $\varphi_3$ have zero mean, the Poincar\'e inequality with $L=\pi R$ yields
\[
    \|\varphi_j\|
    \le
    \frac{L}{\pi}\|\varphi_{j,s}\|
    =
    R\|\varphi_{j,s}\|,
    \quad
    j=1,2,3.
\]
In addition, since $\bv{\varphi}_s\in H^1_0(I)^3$, 
\[
    \|\bv{\varphi}_s\|
    \le
    \frac{L}{\pi}\|\bv{\varphi}_{ss}\|
    =
    R\|\bv{\varphi}_{ss}\|.
\]
Combining these estimates gives the lower bound. This proves (\ref{item2.lem:=pi.proj}).

(\ref{item3.lem:=pi.proj})
For every $a,b\in\R$, the zero means of $\bv{e}_3\cdot\bv{\varphi}$ and $\bv{e}_1\cdot\bv{\varphi}$ imply
\[
    \|\bv{x}-(\bv{x}^R_{0} + a\bv{e}_3 + b\bv{e}_1)\|_2^2
    =
    \|\bv{\varphi}\|_2^2
    + L(p-a)^2
    + L(q-b)^2.
\]
Minimization in $a,b$ proves (\ref{item3.lem:=pi.proj}).

(\ref{item4.lem:=pi.proj})
This follows from \eqref{eq.L.arc.inv}, (\ref{item1.lem:=pi.proj}) above, and Lemmas \ref{lem:L.arc.reduct} and \ref{lem:=pi.inv}. \end{proof}

Now we prove Proposition \ref{prop:=pi.equiv} as a corollary of Lemmas \ref{lem:=pi.coercive} and \ref{lem:=pi.proj}.

\begin{proofx}{Proposition \ref{prop:=pi.equiv}}
Let $\bv{x}\in Y_\pi$ and define $\bv{\varphi}\in X
_\pi$ by \eqref{eq1.lem:=pi.proj}. By Lemma \ref{lem:=pi.proj} (\ref{item4.lem:=pi.proj}), we have $\mcL[\bv{x}]=\mcL[\bv{\varphi}]$. In view of Lemma \ref{lem:=pi.proj} (\ref{item2.lem:=pi.proj}) and (\ref{item3.lem:=pi.proj}), it suffices to show that 
\[
    c\|\bv{\varphi}_{ss}\|^2
    \le
    \mcL[\bv{x}]
    \le
    C\|\bv{\varphi}\|_2^2, 
\]
for some constants $c,C$ depending on $R$. The upper bound follows from the definition of $\mcL$ in \eqref{def.L.arc}. The lower bound follows from Lemma \ref{lem:=pi.coercive}. The proof is complete.
\end{proofx}

    \subsection{Unstable Arc: The Case \texorpdfstring{$\theta>\pi$}{theta > pi}}
    \label{sec.arc.instab}

Let $\pi<\theta<2\pi$. As in Section \ref{sec.<pi.stab}, we define
\[
    \mathsf{A}_{\theta}
    :=
    \bv{x}^R_{0} + \R \bv{e}_3,
    \qquad
    \opdist(\bv{x},\mathsf{A}_{\theta})
    :=
    \inf_{\bv{y}\in\mathsf{A}_{\theta}}
    \|\bv{x}-\bv{y}\|_2.
\]
The aim of this section is to prove the following instability result.

\begin{theorem}
\label{thm:>pi.instab}
There exist constants $\eta,C_0>0$, depending on $R,\theta$, such that the following holds. For every $0<\eps<1/2$, there exist initial data $\bv{x}^{(\eps)}_{0}\in M_\theta\cap C^\infty(\overline{I})^3$ satisfying \eqref{cond1.thm:eq.arc} and
\begin{equation}
\label{est1.thm:>pi.instab}
    \|\bv{x}^{(\eps)}_{0}-\bv{x}^R_{0}\|_{4}
    \le
    C_0 \eps, 
\end{equation}
such that the corresponding solution $\bv{x}^{(\eps)}$ in Theorem \ref{thm:eq.arc} satisfies, at the time
\[
    T^\eps
    :=
    \frac{LR}{\pi}
    \left(
    1-\frac{\pi^2}{\theta^2}
    \right)^{-1/2}
    \log
    \frac{1}{\eps}
    >0,
\]
the estimate
\begin{equation}
\label{est2.thm:>pi.instab}
    \left\|\,
    |\bv{x}_{ss}^{(\eps)}(T^\eps)|
    -
    \frac1R
    \,\right\|
    \ge
    \eta.
\end{equation}
In particular,
\[
    \opdist(\bv{x}^{(\eps)}(T^\eps),\mathsf{A}_{\theta})
    \ge
    \eta.
\]
\end{theorem}

\begin{remark}
\label{rem1:thm:>pi.instab}
An alternative proof of Theorem \ref{thm:>pi.instab} can be based on an analysis of the linearized nonlinear Schr\"odinger (NLS) equation. Under the generalized Hasimoto transformation, the circular arc $\bv{x}^R_{0}$ corresponds, up to a constant phase, to the plane wave
\[
    -\frac1R e^{-\ii t/(2R^2)}
\]
of the focusing cubic NLS equation with Neumann boundary conditions. Writing
\[
    q(s,t)
    =
    e^{-\ii t/(2R^2)}
    \left(
    -\frac1R
    +
    u(s,t)
    \right),
\]
we obtain the linearized equation
\[
    \ii u_t
    =
    u_{ss}
    + \frac1{R^2}\Re u,
    \qquad 
    u_s(0,t)=u_s(L,t)=0.
\]
This equation has an eigenvalue associated with the first mode $\cos(\pi s/L)$, namely,
\[
    \frac{\pi}{LR}
    \left(1-\frac{\pi^2}{\theta^2}\right)^{1/2}>0.
\]
Hence a fixed-point theorem applied to the NLS equation gives a solution approaching the plane wave exponentially as $t\to-\infty$. Then the inverse generalized Hasimoto transformation gives solutions satisfying the conclusions of Theorem \ref{thm:>pi.instab}. Below, we instead use an explicit Akhmediev breather and avoid this perturbation argument.
\end{remark}

\begin{remark}
\label{rem2:thm:>pi.instab}
For the periodic focusing cubic NLS equation, Gustafson, Le Coz, and Tsai \cite[Proposition 3.2 and Corollary 3.9]{GustafsonLeCozTsai2017} gave a variational characterization and proved orbital stability of spatially homogeneous waves. After even extension of the Neumann problem to period $2L$, and in the normalization used here, their stability condition for amplitude $1/R$ becomes $L\le\pi R$, consistent with our stability range in Sections \ref{sec.<pi.stab} and \ref{sec.=pi.stab}.
\end{remark}

We prove Theorem \ref{thm:>pi.instab} by constructing $\bv{x}^{(\eps)}$ from an Akhmediev breather \cite{AkhmedievKorneev1986,AlejoFanelliMunoz2019}. Set
\[
    \mu
    :=
    \frac{\pi}{L},
    \qquad
    \lambda
    :=
    \frac{\pi}{LR}
    \left(
    1-\frac{\pi^2}{\theta^2}
    \right)^{1/2}, 
\]
and
\[
    k
    :=
    \frac1R,
    \qquad
    \rho
    :=
    \frac{\mu}{k}
    =
    \frac\pi\theta,
    \qquad
    b
    :=
    \sqrt{1-\rho^2}.
\]
The Akhmediev breather is expressed as follows: 
\begin{equation}
\label{def.Q}
    Q(s,t)
    :=
    -k e^{-\ii k^2t/2}
    \left(
    1
    -
    \frac
    {2\rho^2\cosh(\lambda t)-2\ii\rho b\sinh(\lambda t)}
    {\cosh(\lambda t)-b\cos(\mu s)}
    \right).
\end{equation}
Note that $0<\rho,b<1$ and $\lambda=k^2\rho b$.

More precisely, if $A(t,x)$ denotes the solution in \cite[(1.1)]{AlejoFanelliMunoz2019} with $a=b^2/2$, then 
\[
    Q(s,t)
    =
    -k
    \overline{
    A\left(\frac{k^2t}{2}, \frac{ks}{\sqrt2}\right)
    }. 
\]
We omit the proof of the following lemma, which follows from the representation \eqref{def.Q}.

\begin{lemma}
\label{lem:Q}
The following hold.
\begin{enumerate}[(i)]
\item
\label{item1.lem:Q}
$Q$ is smooth on $\R\times\R$ and its restriction satisfies 
\[
    \left\{
    \begin{array}{ll}
    \displaystyle{
    \ii Q_t
    =
    Q_{ss}
    +
    \frac12|Q|^2
    Q,
    }
    &(s,t)\in I\times \R, \\[3mm]
    Q_{s}(0,t)
    =
    0,
    \quad
    Q_{s}(L,t)
    =
    0, 
    &t\in \R. 
    \end{array}\right.
\]

\item
\label{item2.lem:Q}
There exists a complex number $C_{-}$ such that $|C_{-}|=1$ and
\[
    \|e^{\ii k^2t/2}Q(\cdot,t) + kC_{-}\|_{2}
    \le
    Ce^{\lambda t},
    \quad
    t\le0.
\]
\end{enumerate}
\end{lemma}

To obtain solutions of \eqref{eq.arc} from the Akhmediev breather $Q$, we employ the inverse generalized Hasimoto transformation; see \cite[Sections 3.1--3.2]{Aiki2023}. Define
\[
    \Psi(s,t)
    :=
    e^{\ii\gamma(t)}
    Q(s,t),    
    \qquad
    \gamma(t)
    :=
    \frac12
    \int_0^t
    |Q(0,\tau)|^2
    \dd\tau.
\]
Then $\Psi$ satisfies
\[
    \left\{
    \begin{array}{ll}
    \displaystyle{
    \ii \Psi_t
    =
    \Psi_{ss}
    +
    \frac12(|\Psi|^2 - |\Psi(0,t)|^2)
    \Psi,
    }
    &(s,t)\in I\times \R, \\[3mm]
    \Psi_{s}(0,t)
    =
    0,
    \quad
    \Psi_{s}(L,t)
    =
    0, 
    &t\in \R. 
    \end{array}\right.
\]
Setting $\Psi:=\Psi_1+\ii \Psi_2$, we define $\bv{v}$, $\bv{e}$, and $\bv{w}$ as the solution of the following system:
\begin{equation}
\label{eq.vew}
    \left\{
    \begin{array}{ll}
    \bv{v}_s
    =
    \Psi_1\bv{e}
    + \Psi_2\bv{w},
    &(s,t)\in I\times \R, \\[3mm]
    \bv{e}_s
    =
    -\Psi_1\bv{v},
    &(s,t)\in I\times \R, \\[3mm]
    \bv{w}_s
    =
    -\Psi_2\bv{v},
    &(s,t)\in I\times \R, \\[3mm]
    (\bv{v},\bv{e},\bv{w})(0,t)
    =
    (\bv{e}_2,\bv{e}_1,\bv{e}_3),
    &t\in \R.
    \end{array}\right.    
\end{equation}
Then $\{\bv{v},\bv{e},\bv{w}\}$ is an orthonormal basis satisfying
\[
    \bv{w}
    =
    \bv{e}\times\bv{v},
    \qquad
    \bv{e}
    =
    \bv{v}\times\bv{w}.
\]
Moreover, $\{\bv{v},\bv{e},\bv{w}\}$ solves 
\[
    \left\{
    \begin{array}{ll}
    \bv{v}_t
    =
    \Psi_{2,s}\bv{e} - \Psi_{1,s}\bv{w}, 
    &(s,t)\in I\times \R, \\[3mm]
    \displaystyle{
    \bv{e}_t
    =
    -\Psi_{2,s}\bv{v} 
    + \frac12(|\Psi|^2 - |\Psi(0,t)|^2)
    \bv{w},
    }
    &(s,t)\in I\times \R, \\[3mm]
    \displaystyle{
    \bv{w}_t
    =
    \Psi_{1,s}\bv{v}
    -
    \frac12(|\Psi|^2 - |\Psi(0,t)|^2)
    \bv{e}
    }
    &(s,t)\in I\times \R. 
    \end{array}\right.
\]
Importantly, $\bv{v}$ satisfies the following.

\begin{lemma}
\label{lem:>pi.v}
The vector $\bv{v}$ above satisfies $|\bv{v}|=1$, $|\bv{v}_s|=|Q|$, 
\[
    \bv{v}(s,0)
    =
    -\sin\alpha(s)\bv{e}_1
    + \cos\alpha(s)\bv{e}_2, 
    \qquad
    \alpha(s)
    :=
    -\int_0^s 
    Q(r,0)
    \dd r,
\]
and
\[
    \left\{
    \begin{array}{ll}
    \bv{v}_{t}
    =
    \bv{v}\times \bv{v}_{ss},
    &(s,t)\in I\times \R, \\[3mm]
    \bv{v}(0,t)
    =
    \bv{e}_{2},
    \quad
    \bv{v}(L,t)
    =
    \bv{b}_\theta, 
    &t\in \R. 
    \end{array}\right. 
\]
Moreover, there exists a constant $C>0$, depending on $R,\theta$, such that
\begin{equation}
\label{est1.lem:>pi.v}
    \|\bv{v}(\cdot,t) - \bv{x}^R_{0,s}\|_{3}
    \le
    C e^{\lambda t},
    \quad
    t\le0.
\end{equation}
\end{lemma}

\begin{proof}
The argument in \cite[Section 3.2]{Aiki2023} applied to $\{\bv{v},\bv{e},\bv{w}\}$ verifies all the assertions except for the condition $\bv{v}(L,t)=\bv{b}_\theta$ and the estimate \eqref{est1.lem:>pi.v}.

For the condition $\bv{v}(L,t)=\bv{b}_\theta$, thanks to \cite[Section 3.2]{Aiki2023} again, we first observe that $\bv{v}(L,t)=\bv{v}(L,0)$ for every $t$. By the definition \eqref{def.Q}, we see that 
\[
    \alpha(L)
    =
    -\int_0^L
    Q(r,0)
    \dd r
    =
    k
    \int_0^L
    \left(
    1
    -
    \frac
    {2\rho^2}
    {1-b\cos(\mu r)}
    \right)
    \dd r.
\]
Using $\mu L=\pi$ and $0<b<1$, we compute 
\[
\begin{split}
    \int_0^L
    \frac{1}{1-b\cos(\mu r)} 
    \dd r
    &=
    \frac{1}{\mu} 
    \int_{0}^{\pi} 
    \frac{1}{1-b\cos s}
    \dd s \\[3mm]
    &=
    \frac{1}{\mu} 
    \lim_{\sigma \nearrow \pi} 
    \bigg[
    \frac{2}{\sqrt{1-b^2}}
    \arctan
    \bigg(
    \sqrt{\frac{1+b}{1-b}}
    \tan\Big(\frac{s}{2}\Big)
    \bigg)
    \bigg]_{0}^{\sigma} \\[3mm]
    &=
    \frac{\pi}{\mu\sqrt{1-b^2}}.
\end{split}
\]
Thus
\[
    \alpha(L)
    =
    kL
    - \frac{2\pi k \rho^2}{\mu\sqrt{1-b^2}} 
    =
    kL
    - \frac{2\pi k\rho}{\mu} 
    =
    \theta - 2\pi, 
\]
which proves, by the definition of $\bv{b}_\theta$ in \eqref{def.b}, 
\[
    \bv{v}(L,0)
    =
    -\sin(\theta-2\pi) \bv{e}_1
    + \cos(\theta-2\pi) \bv{e}_2
    =
    \bv{b}_\theta. 
\]

For the estimate \eqref{est1.lem:>pi.v}, Lemma \ref{lem:Q} (\ref{item2.lem:Q}) and the trace theorem imply that
\[
    \big|
    |Q(0,t)|^2 - k^2
    \big|
    =
    \big|
    |e^{\ii k^2t/2}Q(0,t)|^2 - |kC_{-}|^2 
    \big|
    \le
    Ce^{\lambda t},
    \quad
    t\le0.
\]
Hence $\gamma(t)-k^2t/2$ converges to a real number $\gamma_{-}$ with error $O(e^{\lambda t})$. Setting 
\[
    c_{-}
    :=
    -ke^{\ii\gamma_{-}} C_{-}, 
\]
we have $|c_{-}|=k$ and
\[
    \|\Psi(\cdot,t) - c_{-}\|_{2}
    \le
    C e^{\lambda t},
    \quad
    t\le0.
\]
The system \eqref{eq.vew} with $\Psi$ replaced by $c_{-}$ has the solution
\[
    \bv{v}_{-}(s)
    :=
    \cos(ks)\bv{e}_2
    +
    \frac{\sin(ks)}k
    \big(
    (\Re c_{-})\bv{e}_1
    + (\Im c_{-})\bv{e}_3
    \big).
\]
By differentiating and applying Gronwall's inequality, we obtain
\[
    \|\bv{v}(\cdot,t) - \bv{v}_{-}\|_{3}
    \le 
    C\|\Psi(\cdot,t) - c_{-}\|_{2}
    \le 
    C e^{\lambda t},
    \quad
    t\le0.
\]
Thus the proof of \eqref{est1.lem:>pi.v} is complete if $\bv{v}_{-}=\bv{x}^R_{0,s}$ is shown. Using the trace theorem, 
\[
    \bv{v}_{-}(L)
    =
    \lim_{t\to-\infty}
    \bv{v}(L,t)
    =
    \bv{b}_\theta.
\]
Then, since $kL=\theta$ and $\sin\theta\ne0$, 
\[
    -\sin\theta
    =
    \bv{v}_{-}(L)\cdot\bv{e}_1
    =
    \frac{\sin\theta}k
    \Re c_{-}, 
    \qquad
    0
    =
    \bv{v}_{-}(L)\cdot\bv{e}_3
    =
    \frac{\sin\theta}k
    \Im c_{-}, 
\]
which yields $\Re c_{-}=-k$ and $\Im c_{-}=0$. Therefore $\bv{v}_{-}=\bv{x}^R_{0,s}$ follows. 
\end{proof}

Now we prove Theorem \ref{thm:>pi.instab}.

\begin{proofx}{Theorem \ref{thm:>pi.instab}}
Let $0<\eps<1/2$ and $\bv{v}$ be the vector in Lemma \ref{lem:>pi.v}. Set
\[
    T^\eps
    :=
    \frac1\lambda\log\frac1\eps,
    \qquad
    \bv{h}^\eps(s)
    :=
    \bv{v}(s,-T^\eps) 
    - 
    \bv{x}_{0,s}^R(s).
\]

The initial data $\bv{x}^{(\eps)}_{0}$ are given by 
\[
    \bv{x}^{(\eps)}_{0}(s)
    :=
    \bv{x}^R_{0}(s)
    + 
    d^\eps\bv{e}_1
    + 
    \int_0^s
    \bv{h}^\eps(r)
    \dd r, 
    \qquad
    d^\eps
    :=
    \frac{1}{\sin\theta}
    \int_I
    \bv{h}^\eps\cdot\bv{b}_\theta
    \dd s. 
\]
By construction, 
\[
    |\bv{x}^{(\eps)}_{0,s}(s)|
    =
    1,
    \qquad
    (\bv{x}^{(\eps)}_{0}-\bv{x}^R_{0})(0)\cdot\bv{e}_2
    =
    0,
    \qquad
    (\bv{x}^{(\eps)}_{0}-\bv{x}^R_{0})(L)\cdot\bv{b}_\theta
    =
    0.
\]
Thus $\bv{x}^{(\eps)}_{0}\in M_\theta\cap C^\infty(\overline I)^3$. Since $\bv{v}(0,t)$ and $\bv{v}(L,t)$ are constant in $t$ by \cite[Section 3.2]{Aiki2023},
\[
    (\bv{x}^{(\eps)}_{0,s}\times\bv{x}^{(\eps)}_{0,sss})(0)
    =
    (\bv{x}^{(\eps)}_{0,s}\times\bv{x}^{(\eps)}_{0,sss})(L)
    =
    \bv{0}.
\]
Hence the compatibility conditions in \eqref{cond1.thm:eq.arc} are satisfied. Moreover, by \eqref{est1.lem:>pi.v}, 
\[
    \|\bv{x}^{(\eps)}_{0}-\bv{x}^R_{0}\|_{4}
    \le
    C\|\bv{h}^\eps\|_{3}
    \le
    C_0e^{-\lambda T^\eps}=C_0\eps.
\]

The solutions $\bv{x}^{(\eps)}$ are given by 
\[
    \bv{x}^{(\eps)}(s,t)
    :=
    \bv{x}^{(\eps)}_{0}(s)
    +
    \int_0^t
    (\bv{v}\times\bv{v}_s)(s,\tau-T^\eps)
    \dd\tau,
    \quad
    t\ge0.
\]
Indeed, due to $(\bv{v}\times\bv{v}_s)_s=\bv{v}_t$, we see that \[
    \bv{x}^{(\eps)}_s(s,t)
    =
    \bv{v}(s,t-T^\eps),
    \qquad
    \bv{x}^{(\eps)}_t
    =
    \bv{x}^{(\eps)}_s\times\bv{x}_{ss}^{(\eps)}.
\]
Hence $\bv{x}^{(\eps)}$ is the global-in-time solution of \eqref{eq.arc} in Theorem \ref{thm:eq.arc} by uniqueness. Moreover, using $|\bv{v}_s|=|Q|$ from Lemma \ref{lem:>pi.v} and $k=1/R$, 
\[
    \left\|\,
    |\bv{x}_{ss}^{(\eps)}(\cdot,T^\eps)|
    -
    \frac1R
    \,\right\|
    =
    \left\|\,
    |Q(\cdot,0)|
    -
    k
    \,\right\|
    :=
    \eta
    >
    0, 
\]
where the positivity of $\eta$ follows from $Q(0,0)=k(1+2b)>k$ and the smoothness of $Q$.

In particular, for every $p\in\R$, 
\[
    \|\bv{x}^{(\eps)}(T^\eps)
    -
    (\bv{x}^R_{0}+p\bv{e}_3)\|_2
    \ge
    \|\bv{x}_{ss}^{(\eps)}(T^\eps)
    -
    \bv{x}_{0,ss}^R\| 
    \ge
    \left\|\,
    |\bv{x}_{ss}^{(\eps)}(T^\eps)|
    -
    \frac1R
    \,\right\|
    =
    \eta.
\]
This completes the proof of Theorem \ref{thm:>pi.instab}, with $C_0,\eta$ as above.
\end{proofx}

\begin{remark}
The solutions $\bv{x}^{(\eps)}$ in the proof of Theorem \ref{thm:>pi.instab} satisfy
\[
    \mcL[\bv{x}^{(\eps)}(t)]
    =
    0,
    \quad
    t\ge0,
\]
whereas $\bv{x}^{(\eps)}_{0}\notin\mathsf{A}_{\theta}$ and 
\[
    \|\bv{x}^{(\eps)}_{0} - \bv{x}^R_{0}\|_4
    \to
    0, 
    \quad
    \eps
    \to
    0.
\]
Consequently, when $\pi<\theta<2\pi$, there are no constants $c,r>0$ such that
\[
    \mcL[\bv{x}]
    \ge 
    c\opdist(\bv{x},\mathsf{A}_{\theta})^2,
\]
for every $\bv{x}\in M_\theta$ satisfying $\opdist(\bv{x},\mathsf{A}_{\theta})<r$.

The proof is as follows. By construction, 
\[
    \bv{x}^{(\eps)}_s(s,T^\eps)
    =
    \bv{v}(s,0)
\]
for every $0<\eps<1/2$. Hence, for any $0<\eps,\delta<1/2$, the difference
\[
    \bv{c}(s)
    :=
    \bv{x}^{(\eps)}(s,T^\eps)
    -
    \bv{x}^{(\delta)}(s,T^\delta)
    \in\R\bv{e}_3
\]
is independent of $s$. Moreover,
\[
    \bv{c}
    \cdot
    \bv{e}_2
    =
    \bv{c}
    \cdot
    \bv{b}_\theta=0.
\]
By \eqref{eq.L.arc.inv} and Lemma \ref{lem:L.arc.LOC},
\[
    \mcL[\bv{x}^{(\eps)}(t)]
    =
    \mcL[\bv{x}^{(\eps)}(T^\eps)]
    =
    \mcL[\bv{x}^{(\delta)}(T^\delta)]
    =
    \mcL[\bv{x}^{(\delta)}_{0}],
    \quad
    t\ge0.
\]
Letting $\delta\to 0$, we obtain $\mcL[\bv{x}^{(\eps)}(t)]=\mcL[\bv{x}^R_{0}]=0$ by \eqref{est1.thm:>pi.instab} and the continuity of $\mcL$ on $H^2(I)^3$. If $\bv{x}^{(\eps)}_{0}\in\mathsf{A}_{\theta}$, uniqueness in Theorem \ref{thm:eq.arc} would imply that $\bv{x}^{(\eps)}$ is a translation of $\bv{x}^R_{0}$ along $\bv{e}_3$ for all $t\ge0$. This contradicts the estimate \eqref{est2.thm:>pi.instab}. Hence $\bv{x}^{(\eps)}_{0}\notin\mathsf{A}_{\theta}$.
\end{remark}

    \section{Closed Circles}
    \label{sec.circ}

In this section, we discuss the orbital stability of closed circles. Let $R>0$, $L=2\pi R$, and $\T=\R/L\Z$. We then consider the initial value problem for the LIE:
\begin{align}
\label{eq.circ}
\tag{Circ}
    \left\{
    \begin{array}{ll}
    \bv{x}_{t}
    =
    \bv{x}_{s}\times \bv{x}_{ss},
    &(s,t)\in \T\times (0,\infty), \\[3mm]
    \bv{x}(s,0)
    =
    \bv{x}_{0}(s),
    &s\in \T.
    \end{array}\right. 
\end{align}
The closed circle $\bv{x}^R$ is given by
\[
    \bv{x}^R(s,t)
    :=
    \bv{x}^R_{0}(s)+\frac{t}{R}\bv{e}_3,
    \qquad
    \bv{x}^R_{0}(s)
    :=
    R\cos\left(\frac{s}{R}\right)\bv{e}_1
    +
    R\sin\left(\frac{s}{R}\right)\bv{e}_2.
\]

We recall the solvability of \eqref{eq.circ}. The following result is a consequence of \cite[Theorems 4.3--4.4]{TaniNishiyama1997}.

\begin{theorem}
\label{thm:eq.circ}
Suppose that an initial datum $\bv{x}_0\in H^4(\T)^3$ satisfies
\begin{equation}
\label{cond0.thm:eq.circ}
    |\bv{x}_{0,s}(s)|=1,
    \quad 
    s\in\T.
\end{equation}
Then \eqref{eq.circ} has a unique global-in-time solution $\bv{x}$ in the class
\begin{equation}
\label{est1.thm:eq.circ}
    \bv{x}
    \in
    \bigcap_{j=0}^{2}
    C^j([0,\infty);H^{4-2j}(\T)^3).
\end{equation}
Moreover, $\bv{x}$ satisfies
\begin{equation}
\label{est2.thm:eq.circ}
\begin{split}
    &\bv{x}_s
    \in 
    L^\infty(0,\infty;H^3(\T)^3),\\[3mm]
    &|\bv{x}_s(s,t)|=1,
    \quad 
    (s,t)\in\T\times[0,\infty).
\end{split}
\end{equation}
\end{theorem}

    \subsection{Variational Setting for Closed Circles}
    \label{sec.circ.setting}

We introduce the variational framework for the circle problem \eqref{eq.circ}.

    \subsubsection*{Function Spaces}

Set
\[
    M
    = 
    \{
    \bv{x}\in H^2(\T)^3
    : 
    |\bv{x}_s(s)|=1\ \text{for all }s\in\T
    \}.
\]

    \subsubsection*{Lagrange Functional}

For $\bv{x}\in H^2(\T)^3$, we define
\[
    B[\bv{x}]
    := \frac12\int_{\T}|\bv{x}_{ss}|^2\dd s
\]
and
\[
    \bv{P}[\bv{x}]
    :=
    \frac12\int_{\T}\bv{x}\times\bv{x}_s\dd s,
    \qquad
    P_3[\bv{x}]
    :=
    \bv{P}[\bv{x}]\cdot\bv{e}_3.
\]

Consider the following Lagrange functional $\mcL$ on $H^2(\T)^3$: 
\begin{equation}
\label{def.L.circ}
\begin{split}
    \mcL[\bv{x}]
    &
    :=
    B[\bv{x}]-\frac1{R^3}P_3[\bv{x}] \\
    &=
    \frac12
    \left(
    \int_{\T}|\bv{x}_{ss}|^2\dd s
    -\frac1{R^3}\int_{\T}
    \bv{e}_3\cdot(\bv{x}\times\bv{x}_s)\dd s
    \right).
\end{split}
\end{equation}
The functional $\mcL$ is smooth and quadratic on $H^2(\T)^3$, and its derivative is
\begin{equation}
\label{eq.DL.circ}
    D\mcL(\bv{x})[\bv{\varphi}]
    =
    \int_\T
    \bv{x}_{ss}\cdot\bv{\varphi}_{ss}
    \dd s 
    - \frac1{2R^3}
    \int_\T
    \bv{e}_3\cdot
    (
    \bv{\varphi}\times\bv{x}_s
    + \bv{x}\times\bv{\varphi}_s
    )
    \dd s.
\end{equation}
Periodicity gives $\int_\T\bv{x}_s\dd s=\bv{0}$, and hence
\[
    \bv{P}[\bv{x}+\bv{\tau}]
    =
    \bv{P}[\bv{x}],
    \qquad
    \bv{P}[Q_z^\alpha\bv{x}]
    =
    Q_z^\alpha\bv{P}[\bv{x}], 
\]
where $Q^\alpha_z$ denotes the rotation matrix around the $z$-axis by angle $\alpha$:
\[
    Q_z^\alpha
    = 
    \begin{pmatrix}
        \cos \alpha & -\sin \alpha & 0 \\
        \sin \alpha & \cos \alpha & 0 \\
        0 & 0 & 1
    \end{pmatrix}.
\]
Hence, for every $\bv{x}\in H^2(\T)^3$, $\alpha\in\R$, and $\bv{\tau}\in\R^3$,
\begin{equation}
\label{eq.L.circ.inv}
    \mcL[Q_z^\alpha\bv{x}+\bv{\tau}]
    =
    \mcL[\bv{x}].
\end{equation}
Moreover,
\begin{equation}
\label{eq.L.circ.xR0}
    B[\bv{x}^R_{0}]
    =
    \frac\pi R,
    \qquad
    P_3[\bv{x}^R_{0}]
    =
    \pi R^2,
    \qquad
    \mcL[\bv{x}^R_{0}]
    =
    0.
\end{equation}

The following lemma shows that $\mcL[\bv{x}]$ is conserved for the solutions in Theorem \ref{thm:eq.circ}.

\begin{lemma}
\label{lem:L.circ.LOC}
Let $\bv{x}=\bv{x}(t)$ be the unique global solution of \eqref{eq.circ} in Theorem \ref{thm:eq.circ} with the initial datum $\bv{x}_0\in M \cap H^4(\T)^3$. Then $\bv{x}(t)\in M$ holds for all $t\ge0$. Moreover, $B[\bv{x}(t)]$, $P_3[\bv{x}(t)]$, and $\mcL[\bv{x}(t)]$ are conserved in time.
\end{lemma}

\begin{proof}
The assertion $\bv{x}(t)\in M$ for $t\ge0$ follows from Theorem \ref{thm:eq.circ}. For the conservation of $B[\bv{x}(t)]$, the regularity in \eqref{est1.thm:eq.circ}, periodicity, and $\bv{x}_{st}=\bv{x}_s\times\bv{x}_{sss}$ give
\[
\begin{split}
    \frac{d}{dt}
    B[\bv{x}(t)]
    &=
    \int_\T 
    \bv{x}_{ss}\cdot\bv{x}_{sst}
    \dd s
    =
    -\int_\T 
    \bv{x}_{sss}\cdot\bv{x}_{st}
    \dd s \\
    &=
    -\int_\T 
    \bv{x}_{sss}\cdot
    (\bv{x}_s\times\bv{x}_{sss})
    \dd s
    =
    0.
\end{split}
\]
The conservation of $\bv{P}[\bv{x}(t)]$ is recalled in \cite[Lemma 2.1]{AikiHigaki2026a} and implies that of $P_3[\bv{x}(t)]$. Since $\mcL=B-R^{-3}P_3$, the conservation of $\mcL[\bv{x}(t)]$ follows as well.
\end{proof}

The closed circle $\bv{x}^R_{0}$ is a critical point of $\mcL$.

\begin{lemma}
\label{lem:L.circ.deriv}
For any $\bv{\varphi}\in H^2(\T)^3$, we have
\[
    D\mcL(\bv{x}^R_{0})[\bv{\varphi}]
    =
    0.
\]
\end{lemma}

\begin{proof}
Using periodicity and $\bv{x}^R_{0,ssss}=R^{-4}\bv{x}^R_{0}$, we integrate by parts twice to obtain
\[
    \int_\T
    \bv{x}^R_{0,ss}\cdot\bv{\varphi}_{ss}
    \dd s
    =
    \frac1{R^4}
    \int_\T
    \bv{x}^R_{0}\cdot\bv{\varphi}
    \dd s.
\]
Likewise, since $\bv{x}^R_{0,s}\times\bv{e}_3=R^{-1}\bv{x}^R_{0}$,
\[
\begin{split}
    &
    \frac12
    \int_\T
    \bv{e}_3\cdot
    (
    \bv{\varphi}\times\bv{x}^R_{0,s}
    + \bv{x}^R_{0}\times\bv{\varphi}_s
    )
    \dd s\\[3mm]
    &=
    \int_\T
    \bv{e}_3\cdot
    (
    \bv{\varphi}\times\bv{x}^R_{0,s}
    )
    \dd s
    =
    \frac1R
    \int_\T
    \bv{x}^R_{0}\cdot\bv{\varphi}
    \dd s.
\end{split}
\]
Substituting these identities into the formula \eqref{eq.DL.circ} finishes the proof. 
\end{proof}

Then, by the same argument as in the proof of Lemma \ref{lem:L.arc.reduct}, we have the following.

\begin{lemma}
\label{lem:L.circ.reduct}
For any $\bv{x}=\bv{x}^R_{0} + \bv{\varphi}$ where $\bv{\varphi}\in H^2(\T)^3$, we have
\[
    \mcL[\bv{x}] 
    = 
    \mcL[\bv{\varphi}].
\]
\end{lemma}

    \subsection{Stability of Closed Circles}
    \label{sec.circ.stab}

Define
\[   
    \mathsf{C}
    := 
    \{
    Q_z^\alpha\bv{x}^R_{0} + \bv{\tau}
    :
    \alpha\in\R,\ \bv{\tau}\in\R^3
    \},
    \qquad
    \opdist(\bv{x},\mathsf{C})
    :=
    \inf_{\bv{y}\in\mathsf{C}}
    \|\bv{x}-\bv{y}\|_2.
\]
The aim of this section is to prove the following stability result.

\begin{theorem}
\label{thm:circ.stab}
Let $\bv{x}=\bv{x}(t)$ be the unique global solution of \eqref{eq.circ} in Theorem \ref{thm:eq.circ} with the initial datum $\bv{x}_0\in M\cap H^4(\T)^3$. Then $\bv{x}(t)$ satisfies 
\[
    \bigg(
    \frac{c}{C}
    \bigg)^{1/2}
    \opdist(\bv{x}_0,\mathsf{C})
    \le
    \opdist(\bv{x}(t),\mathsf{C})
    \le
    \bigg(
    \frac{C}{c}
    \bigg)^{1/2}
    \opdist(\bv{x}_0,\mathsf{C}),
    \quad
    t\ge0,
\]
where $c,C>0$ are the constants in Proposition \ref{prop:circ.equiv} below.
\end{theorem}

The proof of Theorem \ref{thm:circ.stab} is omitted, as it follows an argument similar to that of Theorem \ref{thm:<pi.stab} by combining Lemma \ref{lem:L.circ.LOC} with the following proposition.

\begin{proposition}
\label{prop:circ.equiv}
There exist constants $c,C>0$ depending on $R$ such that, for any $\bv{x}\in M$, 
\[
    c\opdist(\bv{x},\mathsf{C})^2
    \le
    \mcL[\bv{x}]
    \le
    C\opdist(\bv{x},\mathsf{C})^2.
\]
In particular, 
\begin{equation}
\label{eq.L.circ.zero}
    \{
    \bv{x}\in M
    :
    \mcL[\bv{x}]=0
    \}
    =
    \mathsf{C}.
\end{equation}
\end{proposition}

\begin{remark}
\label{rem:prop:circ.equiv}
The proof below gives the choices
\[
    c
    =
    \frac1{4(1+R^2+R^4)}, 
    \qquad
    C
    =
    1.
\]
These constants are not asserted to be optimal. The restriction $\bv{x}\in M$ in \eqref{eq.L.circ.zero} is essential: for $\rho\neq\pm1$, the curve $\rho \bv{x}^R_{0}$ satisfies $\mcL[\rho \bv{x}^R_{0}]=0$ but does not belong to $\mathsf{C}$.
\end{remark}

\begin{remark}
\label{rem:circ.uniqueness}
From the proof of Theorem \ref{thm:circ.stab}, one can verify that for circular initial data, unconditional uniqueness holds in the weak solution class considered in \cite[Theorem 4.3 and Remark 4.1]{NishiyamaTani1994}. More precisely, for any $T>0$, every distributional solution of \eqref{eq.circ} in
\[
    L^2(0,T;H^2(\T)^3)
    \cap H^1(0,T;L^2(\T)^3)
\]
with $\bv{x}_0=\bv{x}^R_0$ coincides with $\bv{x}^R$. Indeed, $|\bv{x}_s|=1$ and the fluid impulse are preserved at this regularity, giving $P_3[\bv{x}(t)]=\pi R^2$. The equality case in \cite[Proposition A.1]{AikiHigaki2026a}, together with the positive sign of $P_3[\bv{x}(t)]$, implies $\bv{x}(t)\in\mathsf{C}$ for almost every $t$. Hence \eqref{eq.circ} gives $\bv{x}_t=R^{-1}\bv{e}_3$, and the initial condition yields $\bv{x}=\bv{x}^R$. We omit further details here.
\end{remark}

To prove Proposition \ref{prop:circ.equiv}, we use Fourier series: for $\bv{v}\in L^2(\T)^3$, write
\[
    \bv{v}(s) 
    = 
    \sum_{n\in\Z} 
    \hat{\bv{v}}_n e^{\ii n s / R}, 
    \qquad
    \hat{\bv{v}}_n
    := 
    \frac{1}{2\pi R}
    \int_{\T} \bv{v}(s) e^{-\ii n s / R} \dd s, 
\]
and
\[
    \Proj{n}\bv{v}(s) 
    = 
    \hat{\bv{v}}_n e^{\ii n s / R}. 
\]
Moreover, for $m\neq0$, 
\[
    \Proj{|n|\ge m}
    =
    \sum_{|k|\ge m} \Proj{k},
    \qquad
    \Proj{|n|\le m}
    =
    \sum_{|k|\le m} \Proj{k}. 
\]

We also use the following key quantity, measuring the modulation from $\mathsf{C}$: 
\begin{equation}
\label{def.a}
    \zeta[\bv{v}]
    :=
    \frac1L
    \int_\T 
    \left(
    v_1(s)+ \ii v_2(s)
    \right) 
    e^{-\ii s/R}
    \dd s.
\end{equation}
Indeed, by direct computation, 
\[
    \zeta[\bv{x}^R_{0}]
    =
    R,
    \qquad
    \zeta[Q_z^\alpha\bv{v} + \bv{\tau}]
    =
    e^{\ii\alpha} \zeta[\bv{v}].
\]
In particular, $\zeta[Q_z^\alpha\bv{x}^R_{0}+\bv{\tau}]=e^{\ii\alpha}R$. More precisely, we have the following.

\begin{lemma}
\label{lem:circ.a}
The following hold.
\begin{enumerate}[(i)]
\item
\label{item1.lem:circ.a}
For $\bv{x}\in M$, we have $|\zeta[\bv{x}]|\le R$. The value $R$ is attained if and only if $\bv{x}\in\mathsf{C}$.

\item
\label{item2.lem:circ.a}
For $\bv{v}\in H^2(\T)^3$ and $\alpha\in\R$,
\[
    \|\bv{v}-Q_z^\alpha\bv{x}^R_{0}\|_2^2
    =
    \|\bv{v}\|_2^2
    + \|\bv{x}^R_{0}\|_2^2 
    -2LR
    \left(
    1 + \frac1{R^2} + \frac1{R^4}
    \right)
    \Re
    (
    e^{-\ii\alpha}\zeta[\bv{v}]
    ).
\]
\end{enumerate}
\end{lemma}

\begin{proof}
(\ref{item1.lem:circ.a})
Since $\bv{x}\in M$, we have $\|\bv{x}_s\|^2=L$. Thus, by the definition, $|\zeta[\bv{x}]|\le R$. Suppose next that $|\zeta[\bv{x}]|=R$ and set $\beta:=\arg \zeta[\bv{x}]$. Integration by parts and $\bv{x}^R_{0,ss}=-R^{-2}\bv{x}^R_{0}$ give
\[
\begin{split}
    \|\bv{x}_s-Q_z^\beta\bv{x}^R_{0,s}\|^2
    =
    2L
    - \frac{2}{R^2}
    \int_\T\bv{x}\cdot Q_z^\beta\bv{x}^R_{0}
    \dd s.
\end{split}
\]
Then, by direct computation and the choice of $\beta$, we have
\[
\begin{split}
    \frac{2}{R^2}
    \int_\T\bv{x}\cdot Q_z^\beta\bv{x}^R_{0}\dd s
    =
    \frac{2L}{R}
    \Re
    (
    e^{-\ii\beta}\zeta[\bv{x}]
    )
    =
    \frac{2L}{R}
    |\zeta[\bv{x}]|
    =
    2L.
\end{split}
\]
Hence $\bv{x}-Q_z^\beta\bv{x}^R_{0}$ is a constant vector, so that $\bv{x}\in\mathsf{C}$. The converse follows from $\zeta[\bv{x}^R_{0}]=R$.

(\ref{item2.lem:circ.a})
Computations similar to those in the proof of (\ref{item1.lem:circ.a}) imply the assertion.
\end{proof}

We also need several auxiliary lemmas.

\begin{lemma}
\label{lem:circ.coercive}
For any $\bv{\varphi}\in H^2(\T)^3$ satisfying 
\[
    \bv{x}^R_{0} + \bv{\varphi}
    \in 
    M,
    \qquad
    \zeta[\bv{x}^R_{0} + \bv{\varphi}]
    \in
    [0,\infty), 
\]
we have
\begin{equation}
\label{est.lem:circ.coercive}
    \mcL[\bv{\varphi}]
    \ge
    \frac14
    \|\bv{\varphi}_{ss}\|^2.
\end{equation}
\end{lemma}

\begin{proof}
We use the high- and low-frequency estimates in \cite{AikiHigaki2026a}. Put $\bv{x}:=\bv{x}^R_{0}+\bv{\varphi}$ and set
\begin{equation}
\label{def.E.D}
    E[\bv{\varphi}]
    :=
    \|\bv{\varphi}_{ss}\|^2
    -
    \frac1{R^2}
    \|\bv{\varphi}_s\|^2
    \ge
    0,
    \qquad
    D[\bv{x}]
    :
    =
    \pi R^2 
    - 
    P_3[\bv{x}]
    \ge
    0.
\end{equation}
Then, by the property of $\bv{x}^R_{0}$ and $\|\bv{x}_s\|^2=L=2\pi R$, 
\[
    E[\bv{\varphi}]
    =
    \|\bv{x}_{ss}\|^2 - \frac1{R^2}\|\bv{x}_s\|^2
    =
    \|\bv{x}_{ss}\|^2 - \frac{2\pi}{R}.
\]
Hence Lemma \ref{lem:L.circ.reduct} gives 
\begin{equation}
\label{est1.prf.lem:circ.coercive}
    \mcL[\bv{\varphi}]
    =
    \mcL[\bv{x}]
    =
    \frac12 E[\bv{\varphi}]
    +\frac1{R^3} D[\bv{x}].
\end{equation}

For the high-frequency modes, the estimate in \cite[Lemma 2.2 and its proof]{AikiHigaki2026a} yields
\begin{equation}
\label{est2.prf.lem:circ.coercive}
    \|(\Proj{|n|\ge 2}\bv{\varphi})_{ss}\|^2
    \le
    \frac43 E[\bv{\varphi}].
\end{equation}

For the low-frequency modes, as in \cite[Section 2.2]{AikiHigaki2026a}, write
\[
    \bv{x} 
    = 
    \sum_{n\in\Z} 
    \hat{\bv{x}}_n e^{\ii n s / R}
\]
with
\[
    \hat{\bv{x}}_1
    =
    \lambda\bv{C}_1
    + \mu\overline{\bv{C}_1}
    + \nu|\bv{C}_1|\bv{e}_3,
    \qquad
    \bv{C}_1
    :=
    \frac{R}{2}
    (\bv{e}_1-\ii\bv{e}_2).
\]
Then the definition of $\zeta[\bv{x}]$ in \eqref{def.a} yields
\[
    \zeta[\bv{x}]
    =
    \frac{R}2(\lambda + \mu)
    + \ii \frac{R}2(-\ii\lambda + \ii\mu)
    =
    \lambda R.
\]
Using the assumption $\zeta[\bv{x}]\in[0,\infty)$ and Lemma \ref{lem:circ.a} (\ref{item1.lem:circ.a}), we see that $\lambda$ is real and 
\begin{equation}
\label{est3.prf.lem:circ.coercive}
    0
    \le
    \lambda
    \le
    1.
\end{equation}
On the other hand, the calculation in \cite[Section 2.2]{AikiHigaki2026a} shows that, for general $\lambda\in \C$,  
\begin{equation}
\label{est4.prf.lem:circ.coercive}
\begin{split}
    \|(\Proj{|n|\le 1}\bv{\varphi})_{ss}\|^2
    =
    \frac{2\pi}{R}
    (
    |1-\lambda|^2
    + |\mu|^2
    + |\nu|^2
    ). 
\end{split}
\end{equation}
Moreover, 
\begin{equation}
\label{est5.prf.lem:circ.coercive}
    |\lambda|^2 + |\mu|^2 + |\nu|^2 
    + 
    \frac{1}{R^2} 
    \sum_{|n|\geq 2}
    n^2 |\hat{\bv{x}}_{n}|^2
    = 
    1, 
\end{equation}
and
\begin{equation}
\label{est6.prf.lem:circ.coercive}
    D[\bv{x}]
    \ge
    \pi R^2\big( 2|\mu|^2 + |\nu|^2\big)
    +
    \frac{\pi}{2}
    \sum_{|n|\geq 2}
    n^2 |\hat{\bv{x}}_{n}|^2.
\end{equation}
Therefore, combining 
\eqref{est3.prf.lem:circ.coercive}--\eqref{est6.prf.lem:circ.coercive}, 
we obtain
\begin{equation}
\label{est7.prf.lem:circ.coercive}
\begin{split}
    \|(\Proj{|n|\le 1}\bv{\varphi})_{ss}\|^2
    &\le
    \frac{2\pi}{R}
    \left(
    1-\lambda^2
    + |\mu|^2
    + |\nu|^2
    \right) \\[3mm]
    &=
    \frac{2\pi}{R}
    \bigg(
    2(|\mu|^2 + |\nu|^2)
    + \frac{1}{R^2} 
    \sum_{|n|\geq 2}
    n^2 |\hat{\bv{x}}_{n}|^2
    \bigg) \\[3mm]
    &\le
    \frac{4}{R^3}
    \bigg(
    \pi R^2
    (|\mu|^2 + |\nu|^2)
    + \frac{\pi}{2} 
    \sum_{|n|\geq 2}
    n^2 |\hat{\bv{x}}_{n}|^2
    \bigg) \\[2mm]
    &\le
    \frac{4}{R^3} D[\bv{x}].
\end{split}
\end{equation}

Finally, from
\eqref{est1.prf.lem:circ.coercive}, \eqref{est2.prf.lem:circ.coercive}, and \eqref{est7.prf.lem:circ.coercive}, 
we see that
\[
    \frac14
    \|\bv{\varphi}_{ss}\|^2
    \le 
    \frac14
    \bigg(
    \frac43 E[\bv{\varphi}]
    + \frac{4}{R^3} D[\bv{x}]
    \bigg)
    \le
    \frac12 E[\bv{\varphi}]
    + \frac{1}{R^3} D[\bv{x}]
    =
    \mcL[\bv{\varphi}].
\]
This completes the proof. 
\end{proof}

\begin{lemma}
\label{lem:circ.proj}
For $\bv{x}\in H^2(\T)^3$, let $\bv{\varphi}$ be defined by
\begin{equation}
\label{eq1.lem:circ.proj}
    \bv{\varphi}
    :=
    Q_z^{-\beta}(\bv{x}-\bv{p})
    -
    \bv{x}^R_{0},
\end{equation}
where, with $\zeta[\bv{x}]$ defined by \eqref{def.a}, 
\begin{equation}
\label{eq2.lem:circ.proj}
    \beta
    :=
    \arg \zeta[\bv{x}], 
    \qquad
    \bv{p}
    :=
    \frac{1}{L}
    \int_\T
    \bv{x}.
\end{equation}
We use the convention that $\beta=0$ if $\zeta[\bv{x}]=0$. Then $\bv{\varphi}$ satisfies the following:
\begin{enumerate}[(i)]
\item
\label{item1.lem:circ.proj}
$\int_\T\bv{\varphi}=\bv{0}$ and $\zeta[\bv{x}^R_{0}+\bv{\varphi}]=|\zeta[\bv{x}]|$. Moreover, if $\bv{x}\in M$, then $\bv{x}^R_{0}+\bv{\varphi}\in M$.

\item
\label{item2.lem:circ.proj}
$c\|\bv{\varphi}\|_2 \le \|\bv{\varphi}_{ss}\|\le \|\bv{\varphi}\|_2$ with the constant $c=(1+R^2+R^4)^{-1/2}$.

\item
\label{item3.lem:circ.proj}
$\opdist(\bv{x},\mathsf{C}) = \|\bv{\varphi}\|_2$.

\item
\label{item4.lem:circ.proj}
$\mcL[\bv{x}] = \mcL[\bv{\varphi}]$.
\end{enumerate}
\end{lemma}

\begin{proof}
(\ref{item1.lem:circ.proj})
Since $\int_\T\bv{x}^R_{0}=\bv{0}$, the definition of $\bv{p}$ gives $\int_\T\bv{\varphi}=\bv{0}$. Moreover,
\[
    \zeta[\bv{x}^R_{0} + \bv{\varphi}]
    =
    \zeta[Q_z^{-\beta}(\bv{x}-\bv{p})]
    =
    e^{-\ii\beta}
    \zeta[\bv{x}]
    =
    |\zeta[\bv{x}]|.
\]
The last assertion follows from
$|(Q_z^{-\beta}(\bv{x}-\bv{p}))_s|=|\bv{x}_s|$.

(\ref{item2.lem:circ.proj})
The proof is similar to that of Lemma \ref{lem:=pi.proj} (\ref{item2.lem:=pi.proj}) and thus omitted.

(\ref{item3.lem:circ.proj})
For every $\alpha\in\R$ and $\bv{\tau}\in\R^3$, the zero means of $\bv{x}-\bv{p}$ and $Q_z^\alpha\bv{x}^R_{0}$ imply
\[
    \|\bv{x}-(Q_z^\alpha\bv{x}^R_{0}+\bv{\tau})\|_2^2
    =
    \|\bv{x}-\bv{p}-Q_z^\alpha\bv{x}^R_{0}\|_2^2
    +
    L|\bv{p}-\bv{\tau}|^2.
\]
Thus the infimum over $\bv{\tau}$ is attained at $\bv{\tau}=\bv{p}$. Lemma \ref{lem:circ.a} (\ref{item2.lem:circ.a}) implies
\[
\begin{split}
    &\|\bv{x}-\bv{p}-Q_z^\alpha\bv{x}^R_{0}\|_2^2 \\
    &=
    \|\bv{x}-\bv{p}\|_2^2
    + \|\bv{x}^R_{0}\|_2^2 
    -2LR
    \left(
    1 + \frac1{R^2} + \frac1{R^4}
    \right)
    \Re\big(
    e^{-\ii\alpha}\zeta[\bv{x}-\bv{p}]
    \big).    
\end{split}
\]
Thus the infimum over $\alpha$ is attained at $\beta$; if $\zeta[\bv{x}]=0$, every $\alpha$ is a minimizer. Hence 
\[
    \opdist(\bv{x},\mathsf{C})
    =
    \|\bv{x}-\bv{p}-Q_z^\beta\bv{x}^R_{0}\|_2
    =
    \|\bv{\varphi}\|_2.
\]

(\ref{item4.lem:circ.proj})
This follows from \eqref{eq.L.circ.inv} and Lemma \ref{lem:L.circ.reduct}.
\end{proof}

Now we prove Proposition \ref{prop:circ.equiv} as a corollary of Lemmas \ref{lem:circ.coercive} and \ref{lem:circ.proj}.

\begin{proofx}{Proposition \ref{prop:circ.equiv}} 
Let $\bv{x}\in M$ and define $\bv{\varphi}$ by \eqref{eq1.lem:circ.proj}. By Lemma \ref{lem:circ.proj} (\ref{item4.lem:circ.proj}), we have $\mcL[\bv{x}]=\mcL[\bv{\varphi}]$. In view of Lemma \ref{lem:circ.proj} (\ref{item2.lem:circ.proj}) and (\ref{item3.lem:circ.proj}), it suffices to show that 
\[
    c\|\bv{\varphi}_{ss}\|^2
    \le
    \mcL[\bv{x}]
    \le
    C\|\bv{\varphi}\|_2^2, 
\]
for some constants $c,C$ depending on $R$. The upper bound follows from the definition of $\mcL$ in \eqref{def.L.circ}. The lower bound follows from Lemma \ref{lem:circ.proj} (\ref{item1.lem:circ.proj}) and Lemma \ref{lem:circ.coercive}. The proof is complete.
\end{proofx}

    \section{Application}
    \label{sec.appl}

In this section, we illustrate how the functional $\mcL$ in Section \ref{sec.circ} remains useful when the conservation laws in Lemma \ref{lem:L.circ.LOC} no longer hold. This is in the spirit of \cite[Section 1A]{GallaySverak2024a}.

We begin with the following lemma, retaining the notation of Section \ref{sec.circ}.

\begin{lemma}
\label{lem:circ.nonneg}
For any $\bv{x}\in H^2(\T)^3$, we have
\begin{equation}
\label{est.lem:circ.nonneg}
    0
    \le
    \mcL[\bv{x}]
    \le
    \|\bv{x}_{ss}\|^2.
\end{equation}
\end{lemma}

\begin{proof}
The Poincar\'e inequality yields
\[
    2|P_3[\bv{x}]|
    \le
    \left\|
    \bv{x}
    -
    \frac{1}{L}
    \int_\T
    \bv{x}
    \dd s
    \right\|
    \,
    \|\bv{x}_s\|
    \le 
    R^3
    \|\bv{x}_{ss}\|^2.
\]
Then
\[
    \frac12
    \|\bv{x}_{ss}\|^2
    -
    \frac{1}{R^3}
    |P_3[\bv{x}]|
    \le
    \mcL[\bv{x}]
    \le
    \frac12
    \|\bv{x}_{ss}\|^2
    + 
    \frac{1}{R^3}
    |P_3[\bv{x}]|
\]
leads to the assertion. 
\end{proof}

We have an a priori estimate for vortex filaments under arclength-preserving forcing:

\begin{proposition}
\label{prop:appl}
Let $T>0$ and 
\[
    \bv{x}
    \in 
    C([0,T];H^4(\T)^3)
    \cap
    C^1([0,T];H^2(\T)^3)
\]
solve
\begin{equation}
\label{eq.prop:appl}
    \bv{x}_t
    =
    \bv{x}_s\times\bv{x}_{ss}
    + \bv{f},
    \qquad
    \bv{x}(0)
    =
    \bv{x}_0\in M,    
\end{equation}
where $\bv{f}\in C([0,T];H^2(\T)^3)$. Suppose that
\[
    \bv{x}_s\cdot\bv{f}_s
    =
    0
\]
holds on $\T\times[0,T]$. Then $\bv{x}(t)\in M$ for all $t\in[0,T]$ and
\[
    \opdist(\bv{x}(t),\mathsf{C})
    \le
    \sqrt{\frac{C}{c}}
    \opdist(\bv{x}_0,\mathsf{C})
    +
    \frac{1}{\sqrt{c}}
    \int_0^t
    \|\bv{f}_{ss}(\tau)
    \|\dd\tau.
\]
Here $c,C>0$ are the constants in
Proposition~\ref{prop:circ.equiv}.
\end{proposition}

\begin{proof}
Differentiating the equation \eqref{eq.prop:appl} with respect to $s$ gives
\[
    \frac12\partial_t|\bv{x}_s|^2
    =
    \bv{x}_s\cdot\bv{f}_s
    =
    0,
\]
so that $\bv{x}(t)\in M$. Integration by parts and $|\bv{x}_s|=1$ give
\[
    \frac{d}{dt}
    B[\bv{x}]
    =
    \int_\T
    \bv{x}_{ss}\cdot\bv{f}_{ss}
    \dd s,
    \qquad
    \frac{d}{dt}
    \bv{P}[\bv{x}]
    =
    \int_\T
    \bv{f}\times\bv{x}_s
    \dd s.
\]
Consequently, using the notation in \eqref{eq.DL.circ}, 
\[
    \frac{d}{dt}
    \mcL[\bv{x}(t)]
    =
    D\mcL(\bv{x}(t))
    [\bv{f}(t)].
\]
Since $\mcL$ is nonnegative by Lemma \ref{lem:circ.nonneg}, we have
\[
    |D\mcL(\bv{x}(t))[\bv{f}(t)]|
    \le
    2
    \sqrt{\mcL[\bv{x}(t)]}
    \,
    \sqrt{\mcL[\bv{f}(t)]}. 
\]
Integrating the
inequality for $\sqrt{\mcL[\bv{x}(t)]+\eps}$, then letting
$\eps\to0$, we obtain
\[
    \sqrt{\mcL[\bv{x}(t)]}
    \le
    \sqrt{\mcL[\bv{x}_0]}
    +
    \int_0^t
    \sqrt{\mcL[\bv{f}(\tau)]}
    \dd\tau.
\]
Combining Lemma \ref{lem:circ.nonneg} and Proposition \ref{prop:circ.equiv} proves the desired estimate.
\end{proof}

The condition $\bv{x}_s\cdot\bv{f}_s=0$ does not hold for general prescribed external flows $\bv{f}(s,t)=\bv{F}(\bv{x}(s,t),t)$. Without this condition, an appropriate model introduced in \cite{AikiIguchi2019} is
\[
    \bv{x}_t
    =
    \frac
    {\bv{x}_\xi\times\bv{x}_{\xi\xi}}
    {|\bv{x}_\xi|^3}
    +
    \bv{F}(\bv{x},t).
\]
Local well-posedness was proved in \cite[Theorem 2.1]{AikiIguchi2019} under suitable assumptions. Extending our results to that setting requires control of stretching and is not pursued here.

\end{document}